\documentclass[11pt]{article}

\usepackage{amsthm,amsmath,amssymb,amsfonts}
\usepackage{color}
\usepackage{a4wide}
\usepackage{graphicx}
\newtheorem{tw}{Theorem}
\newtheorem{prop}{Proposition}
\newtheorem{lem}{Lemma}
\newtheorem{df}{Definition}
\newtheorem{fa}{Fact}
\newtheorem{rem}{Remark}

\newcommand{\vt}{\vartheta}
\newcommand{\ux}{\underline{x_1}}
\newcommand{\uxd}{\underline{x_2}}
\newcommand{\uy}{\underline{y_1}}
\newcommand{\uxx}{\underline{x}}
\newcommand{\uyy}{y^{in}}
\newcommand{\ox}{\overline{x_1}}

\newcommand{\de}{\partial}
\newcommand{\ep}{\epsilon}
\newcommand{\dx}{{\rm dx}}
\renewcommand{\qed}{\begin{flushright} $\square$ \end{flushright}}

\renewcommand{\div}{{\rm div\,}}
\renewcommand{\d}{\partial}
\newcommand{\rot}{{\rm rot\,}}
\newcommand{\R}{\mathbb{R}}

\begin{document}
\begin{center}
{\Large Stationary compressible Navier - Stokes Equations with inflow condition in domains with piecewise analytical boundaries}

\vskip10mm

{\bf Piotr B. Mucha$^1$, \quad Tomasz Piasecki$^2$}

\vskip3mm

\emph{Institute of Applied Mathematics and Mechanics, University of Warsaw}

\vskip2mm

\emph{Banacha 2, 02-097 Warsaw}

\vskip3mm

$^1$ e-mail: p.mucha@mimuw.edu.pl
$^2$ e-mail: tpiasecki@mimuw.edu.pl

%\vskip5mm

%{\bf Tomasz Piasecki}

%\vskip3mm

%\emph{Institute of Mathematics, Polish Academy of Sciences}

\vskip10mm

\end{center}



\vskip10mm

\centerline{\bf Abstract}

We show the existence of strong solutions in Sobolev-Slobodetskii spaces to the stationary compressible 
Navier-Stokes equations
with inflow boundary condition. Our result holds provided certain condition on the shape of the boundary around 
the points where characteristics of the continuity equation are tangent to the boundary, which
holds in particular for piecewise analytical boundaries. The mentioned situation creates a singularity which limits regularity
at such points. 
We show the existence and uniqueness of regular solutions in a vicinity of given laminar solutions 
under the assumption that the pressure is a linear function of the density. 
The proofs require the language of suitable fractional Sobolev spaces.
In other words our result is an example where application of fractional spaces is irreplaceable, 
although the subject is a classical system.

\vskip10mm

\section{Introduction}
We investigate the existence of regular solutions to stationary barotropic compressible Navier-Stokes
equations in a two dimensional bounded domain $\Omega$ with nonzero inflow/outflow through the  
boundary. The complete system reads
\begin{equation} \label{main_system}
\begin{array}{lcr}
\rho v \cdot \nabla v -\mu \Delta v - (\mu+\nu) \nabla {\rm div}\; v
+\nabla \pi(\rho) =0 & \mbox{in} & \Omega,\\
{\rm div}\;(\rho v)=0 & \mbox{in} & \Omega,\\
n\cdot 2 \mu \mathbf{D}(v) \cdot \tau +f v\cdot \tau=b, &
\mbox{on} &\Gamma,\\
n \cdot v=d & \mbox{on} & \Gamma,\\
\rho=\rho_{in} & \mbox{on} & \Gamma_{in},
\end{array}
\end{equation}
where the velocity field of the fluid $v$ and the density $\rho$ are the unknown functions describing the flow.
We distinguish the parts of the boundary:
\begin{equation} \label{def_bdry}
\begin{array}{c}
\Gamma_{in} = \{x \in \Gamma: d < 0 \}, \quad 
\Gamma_{out} = \{x \in \Gamma: d > 0 \},\\[3pt] 
\Gamma_{0} = \{x \in \Gamma: d = 0 \}, \quad
\Gamma_* = \overline{\Gamma_0} \cap \overline{\Gamma_{in} \cup \Gamma_{out}}. 
\end{array}
\end{equation}  
We show the existence of a solution in fractional Sobolev spaces $u \in W^{1+s}_p(\Omega)$, $\rho \in W^s_p(\Omega)$,
where $u$ is the velocity field of the fluid and $\rho$ is the density.
Our choice of functional spaces allows to overcome the problem of singularity in the continuity equation
and obtain boundedness of the density. Before we formulate the problem more precisely we give a brief
overview of the state of art in the topic, focusing on the scope of interest of this paper, that is
on regular stationary solutions, mentioning also the most important results concerning global 
weak stationary solutions.
For more complete overview of known results in the mathematical theory of compressible flows we refer
to the monographs \cite{NoS} and \cite{PSbook}.

The mathematical theory of stationary solutions to the Navier-Stokes equations describing compressible flows started to develop in early 80's with certain results
on the existence of regular solutions, first in Hilbert spaces (\cite{Valli1}) and later in $L_p$ framework (\cite{BdV}).
However, all of these results required certain smallness assumptions on the data and concerned
mostly homogeneous boundary conditions with vanishing normal component of the velocity.

In the 90's the famous
result of Lions \cite{Li} on the existence of weak solutions for homogeneous Dirichlet boundary conditions 
triggered the development of global existence theory of weak solutions.
The result was improved by Novo and Novotn\'y \cite{NoNo} who adopted the nonsteady  
approach of Feireisl et al. \cite{FeNoPe} and then extended by Mucha and Pokorn\'y  
to the case of slip boundary conditions in the barotropic case (\cite{PMMP1},\cite{PMMP2})
and for the full system including thermal effects \cite{PMMP3}, see also the result 
for a system involving radiation effects in \cite{KNP}.
Further improvements in the theory of regular solutions has been made
in the nineties (\cite{NoPa},\cite{NoPi}) but mostly for homogeneous boundary data.

It should be emphasized that all above mentioned global results concern the case of normal component
of the velocity vanishing on the boundary. 
If the normal component of the velocity does not vanish, 
substantial mathematical difficulties arises in the analysis of the continuity equation, which 
can be reduced to a stationary transport equation. 
Namely, the hyperbolicity of the continuity equation makes it necessary to prescribe the density on the part of the boundary where the fluid enters the domain,
called briefly the inflow part. 
Solvability of either a time-dependent 
or stationary transport equation is of utmost importance in the mathematical analysis of the Navier-Stokes
equations for compressible flow. For a recent application of the theory of transport equation in the context 
of the existence of weak solutions to the compressible Navier-Stokes equation we can refer to \cite{BJ}. 
Important developments in the theory of transport equation, strenghtening the classical results 
of DiPerna and Lions \cite{DPL}, has been made in \cite{CDL} and \cite{M1}.  
 
The mathematical investigation of inflow/outflow problems began with the work of Valli and Zaj\c aczkowski
\cite{VZ},
who investigated time-dependent problem obtaining also an existence result in the stationary case.
Then the development of existence theory for inhomogeneous boundary data has has been hindered
by mathematical difficulties on the one hand and the interest turned mostly towards global existence 
of weak solutions on the other, until the work by Kweon and Kellogg \cite{Kw1}.
More recently, the existence theory has been developed motivated by applications in shape optimization 
by Plotnikov, Ruban and Sokolowski (\cite{PRS1},\cite{PRS2} and the monograph \cite{PSbook}).
All above results require certain smallness assumptions. Concerning large data problems,
there are only few particular results on global existence of weak solutions for nonstationary problems,
see \cite{Giri}. In the stationary case, due to nontrivial boundary terms it has been for a long time impossible to get basic a priori estimates
and further problems are encountered with the issue of existence and uniqueness for the continuity equation.
The first global existence result has been obtained very recently by Feireisl and Novotn\'y \cite{FN}. 
under the assumption that the pressure is a nondecreasing $C^1$ function of the density satisfying 
${\rm lim}_{\rho \to \bar \rho} \, p(\rho) = +\infty$ for some positive constant $\bar \rho$. 
The proof is based on appropriate regularization of both continuity and momentum equations. 
The key estimates for the approximate systems are obtained using a suitable extension 
of the boundary velocity which is constructed in such a way that it satisfies certain smallness condition even
though the data can be arbitrarily large.   

At first glance a natural functional space for regular solutions is $W^1_p$ for the density and $W^2_p$
for the velocity.  
A regular solution is then understood as a function with weak derivatives satisfying the equations 
almost everywhere. 
However, except some special classes of domains we are not able to obtain the solutions
in the above class for arbitrarily large $p$
(see \cite{Kw1}). 
%and \cite{PS} - note that the limitation on $p$ in both cited papers, although formulated
%in a different way is indeed the same, what is a strong evidence of its optimality)
The reason is a singularity arising in the solution of steady transport equation around the points
where characteristics of this hyperbolic equation become tangent to the boundary,
we refer to these points as singularity points.

On the other hand, the range $p>n$ is important since it gives boundedness of the density due to the
imbedding theorem.
The result from \cite{Kw1} 
%and \cite{PS} 
cover a part of this range, namely $2<p<3$. However,
further increase of $p$ is impossible even under relaxation of the 
boundary singularity. Further investigation of this singularity is therefore an interesting question in view of the development
of the theory of regular solutions.

One possible way to obtain existence for $n<p<\infty$ is to investigate some special domains, such
as a cylindrical domain in \cite{PMTP},\cite{P1},\cite{P2},\cite{Guo} for barotropic case and \cite{PP} for system with thermal
effects or an unbounded domain contained between two parallel planes in \cite{Kw2}.
A possible way to overcome the singularity problem described above in a general domain is an appropriate
choice of functional spaces. In \cite{PRS1}, \cite{PRS2} the existence and uniqueness of solutions
in fractional Sobolev spaces (velocity in $W^{1+s}_p$ and density in $W^s_p$) is shown 
under certain assumption relating the inflow velocity and shape of the boundary around the 
singularity points. However, this result requires additional assumption that the gradient 
of the density and the second gradient of the velocity are in $L_2$. 
 
In this paper we show existence of solutions in fractional Sobolev spaces as above. 
However we do not require the existence of $\nabla \rho$ and $\nabla^2 u$.
Our analysis shows that we are able to show existence of the solutions for $sp>n$ which
gives boundedness of the density. We need to impose a certain limitation on the boundary around the
singularity points, however this assumption is weaker than in \cite{Kw1} and \cite{PRS1},\cite{PRS2} 
and turns out quite natural, in particular it is satisfied by analytical boundaries.

%Obviously, our solution is weaker than a regular one according to the above definition.
%In particular, certain quantities in the equations are no longer given as functions.
%Our solutions lies between weak solutions and regular ones. This is  a price we pay 
%for relaxation of the assumptions on the domain {\bf z tego trzeba sie dobrze wytlumaczyc}

The only result giving uniqueness of solutions to the compressible Navier-Stokes equations for large 
data without information on $\nabla \rho$ is \cite{Hoff} for time-dependent problem. The key 
idea there is to show uniqueness in quite low regularity, namely $L_2$ for the velocity 
and $H^{-1}$ for the density. Then we have to estimate $H^{-1}$ - norm of the pressure 
with $H^{-1}$ norm of the density and for this purpose it is required that the pressure is a linear function 
of the density (or satisfies slightly more general constraint - see (1.16) in \cite{Hoff}).
 
Our result, which is to our knowledge a 
first one giving uniqueness of solutions without any information on the gradient of density in the stationary case,  
combines the ideas from \cite{Hoff} in the context of uniqueness with the approach used 
to obtain existence of regular solutions in the series of papers \cite{P1}, \cite{PMTP} and \cite{PP}. 
It shows that the choice of fractional Sobolev spaces is in a sense natural for considered problem and therefore 
is not only of purely mathematical interest. In particular it may indicate a possible direction 
for the development of the theory of global existence.

\subsection{Functional spaces.} 
In order to formulate more precisely the problem and main result we shall first
recall the definitions of the functional spaces we apply.
We use standard Sobolev spaces $W^k_p$ with natural $k$, which consist 
of functions with weak derivatives up to order $k$ in $L_p(\Omega)$, for the definition we refer
for example to \cite{Ad}. 
However, most important for our result are Sobolev-Slobodetskii spaces $W^s_p$ with fractional $s$. 
For the sake of completeness we recall the definition here (see \cite{Tr}, Section 4.4).
By $W^s_p(\Omega)$ we denote the space of functions for which the norm:
%$$
%\|f\|_{W^s_p} = \left( \int_{\mathbb{R}^2} dx \int_{(0,1)^2} \frac{ f(x+h) - f(x) }{|h|^{2+sp}} \,dh \right)^{1/p}
%$$ 
%is finite ([cite?]). 
%The space $W^s_p$ is usually defined on the whole space. For our application we need a definition on a bounded
%domain of $\mathbb{R}^2$. It can be formulated in the same way, the only difference is that the internal integration 
%can not take us out of the domain. Hence we define
%
\begin{equation} \label{wsp_omega}
\|f\|_{W^s_p(\Omega)} = (\int_\Omega |f|^pdx)^{1/p}+ \left( \iint_{\Omega^2} \frac{ |f(x) - f(y)|^p }{|x-y|^{2+sp}} \,dxdy \right)^{1/p}
\end{equation}
is finite and $s\in (0,1)$. Furthermore, $W^{1+s}_p(\Omega)$ is a space of functions with first weak 
derivatives in $W^s_p(\Omega)$ with the norm
\begin{equation} \label{wsp_omega1}
\|f\|_{W^{1+s}_p(\Omega)} = \|f\|_{W^s_p(\Omega)} + \|\nabla f\|_{W^s_p(\Omega)}. 
\end{equation}
   
%where $\delta=\delta(x)$ depends on the distance of $x$ from the boundary. This dependence is not a problem
%in the above definition since integrability in \eqref{wsp_omega} is determined by the behaviour near $0$.

Let us recall two important features of Sobolev-Slobodetskii spaces. We formulate them in a simplified way
convenient for our applications. 
\begin{fa} \label{imbed}
Assume $\Omega \subset \R^2$ be bounded with $\de \Omega \in C^2$,
Let $u \in W^s_p$ with $sp>n$. Then for $\forall q \leq \infty, \quad \delta>0$  
\begin{align} 
& \|u\|_{L_q} \leq C(q,\Omega) \|u\|_{W^s_p}, \label{sob}\\
& \|u\|_{L_q} \leq \delta \|u\|_{W^s_p} + C(\delta) \|u\|_{L_2}. \label{int} 
\end{align}
\end{fa}
The proofs of more general versions of above facts can be found in \cite{Tr}. 
Furthermore, for given $\ep>0$ let us denote 
$$
\Omega_{\ep}=\{ x \in \Omega: \; {\rm dist}(x,\de \Omega)>\ep\}.
$$
Then by $B^r_{p,\infty}(\Omega)$ we denote a space of functions for which (see \cite{Tr}, Sec. 4.4 or \cite{DevSh}, Sec. 2)
\begin{equation} \label{norm:brp}
\|u\|_{B^r_{p,\infty}(\Omega)}={\rm sup}_{0<|h|\leq h_0}|h|^{-r}\|u(\cdot + h)-u(\cdot)\|_{L_p(\Omega_{|h|})}<\infty,
\end{equation}
where $h_0 = {\rm sup}\{|h|:|\Omega_h|>0\}$.
The for $\Omega$ bounded with $\de \Omega \in C^2$ we have (see \cite{Tr}, Sec. 4.6) 
\begin{equation} \label{imbed:2}
W^r_p(\Omega) \subset B^r_{p,\infty}(\Omega).  
\end{equation} 
Finally, let us denote 
\begin{equation} \label{def_V}
V = \{ v \in W^1_2(\Omega): v \cdot n|_{\Gamma}=0 \}.
\end{equation}

\subsection{Problem formulation}
Let us move to a precise statement of the problem under consideration.
We investigate stationary flow of a barotropic fluid in a 
two dimensional, bounded domain described by the system \eqref{main_system}.
The system is supplied with inhomogeneous slip boundary conditions on the velocity. 
In particular, the normal component of the velocity
does not vanish and, as explained above, we have to prescribe the density on the part of the boundary where the flow
enters the domain. 

Let us have a closer look at the definition of different parts of the boundary \eqref{def_bdry}.
We see that $\Gamma_*$ consist of points where the characteristics of the continuity equation 
\eqref{main_system}$_2$ become tangent to the boundary, we will call it the set of singularity points.
Moreover, let $\Gamma_{in}^s$ be a certain boundary neighborhood of the set of singularity points.
Formally we define it as 
\begin{equation} \label{def_gamma_s}
\Gamma^s_{in} = \{ x\in \Gamma_{in}: {\rm dist}(x,\Gamma_*)<2\eta \},
\end{equation}
where $\eta>0$ is a small number to be precised later.
%Obviously $\Gamma_s$ has nonempty intersection with the parts defined in \eqref{def_bdry}. 
%For simplicity we assume for the moment that $\Gamma_0$
%onsist of two points $x_*,x^*$  which we call the singularity points. Around the singularity points the
%boundary is given as a graph of a function, we denote it by $x_2^l(x_1)$ and $x_2^u(x_1)$. 
%We assume also that the inflow and outflow parts are given by the 
%functions $\Gamma_{in} = \underline{x_1}(x_2)$ and $\Gamma_{out} = \overline{x_1}(x_2)$.  

Our goal is to show the existence of a solution $(u,\rho) \in W^{1+s}_p \times W^s_p$ to the system (\ref{main_system}),
where $s>\frac{n}{p}$, 
which is close to the constant flow 
\begin{equation} \label{const}
(\bar v, \bar \rho) \equiv ([v^*,0], 1),
\end{equation}
where $v^*$ is a positive constant. 
Our method works for a wider class of solutions in which $x_1$ is in a sense dominating direction.
Our motivation for the choice of this fractional order space for the density 
has been explained above;
we want to solve the problem of singularity in the solution of the continuity equation around the singularity points.
On the other hand, by the imbedding theorem we have $W^s_p(\Omega) \in L_\infty(\Omega)$.
Then the choice of the space $W^{1+s}_p$ for the velocity follows naturally from the structure of \eqref{main_system}.  
Obviously such solution no longer satisfies the equations almost everywhere. Therefore
in order to define the solutions we need a weak formulation of the problem \eqref{main_system}.
A natural way is to multiply \eqref{main_system}$_1$ by $\psi \in V$ and \eqref{main_system}$_2$ by a
smooth function $\phi$ such that $\phi|_{\Gamma \setminus \Gamma_{in}=0}$, integrate by parts
and apply boundary conditions \eqref{main_system}$_{3,4}$. 
However, because of inhomogeneous condition \eqref{main_system}$_{4}$ we would obtain boundary terms
with derivatives of $v$. Therefore we first remove the inhomogeneity and introduce weak formulation 
of the perturbed problem \eqref{system}. We obtain the following definition
\begin{df} \label{def_sol} 
A strong solution to the system \eqref{main_system} is a couple
$(v,\rho) \in W^{1+s}_p \times W^s_p$ such that $v = \bar v + u+ u_0$ and $\rho = \bar \rho+w$
where $(u,w)$ is a solution to the system \eqref{system} in the sense of Definition \ref{def_sol1}
and $u_0$ is defined in \eqref{def_ext}.
\end{df}

In order to formulate our main result let us 
introduce the following quantity to measure the distance
of the data of the problem \eqref{main_system} from $(\bar v, \bar \rho)$:
\begin{equation} \label{D0}
D_0 = \|b - f \tau_1\|_{W^{s-1/p}_p(\Gamma)} + \|d - n_1\|_{W^{1+s-1/p}_p(\Gamma)}
+ \|\rho_{in}-1\|_{W^s_p(\Gamma_{in})}.
\end{equation} 
As formulation our main result involves certain properties of the boundary, it is useful
to describe first the domain introducing the necessary assumptions. 

\subsection{The domain. Representative case.} 
Conducting the proof for a general domain with multiple singularity points would lead to unnecessary complications
which would likely hide the main ideas. For clarity of the proof we consider a simple domain with
inflow and outflow parts of the boundary given by 
\begin{equation} \label{dom:1}
\begin{array}{c}
\Gamma_{in}=\{(\underline{x_1}(x_2),x_2): x_2 \in (0,b)\},\\[5pt]
\Gamma_{out}=\{(\overline{x_1}(x_2),x_2): x_2 \in (0,b)\},\\ 
\end{array}
\end{equation}
with 
\begin{equation} \label{def_sing3}
\begin{array}{c}
{\rm lim}_{x_2 \to 0^+} \underline{x_1}(x_2) = {\rm lim}_{x_2 \to b^-} \underline{x_1}(x_2) = 
{\rm lim}_{x_2 \to 0^+} \overline{x_1}(x_2) = {\rm lim}_{x_2 \to b^-} \overline{x_1}(x_2) = 0
\end{array}
\end{equation}
and
\begin{equation} \label{def_sing2}
\begin{array}{c}
{\rm lim}_{x_2 \to 0^+} \underline{x_1}'(x_2) = {\rm lim}_{x_2 \to b^-} \overline{x_1}'(x_2)= -\infty,\\[5pt]
{\rm lim}_{x_2 \to 0^+} \overline{x_1}'(x_2) = {\rm lim}_{x_2 \to b^-} \underline{x_1}'(x_2)= +\infty.
\end{array}
\end{equation}
Then we have two singularity points:
$$
\Gamma_*=\Gamma_0=\{(0,0),(0,b)\}.    
$$
\begin{figure} 
\begin{center}
\includegraphics[width=0.7\textwidth]{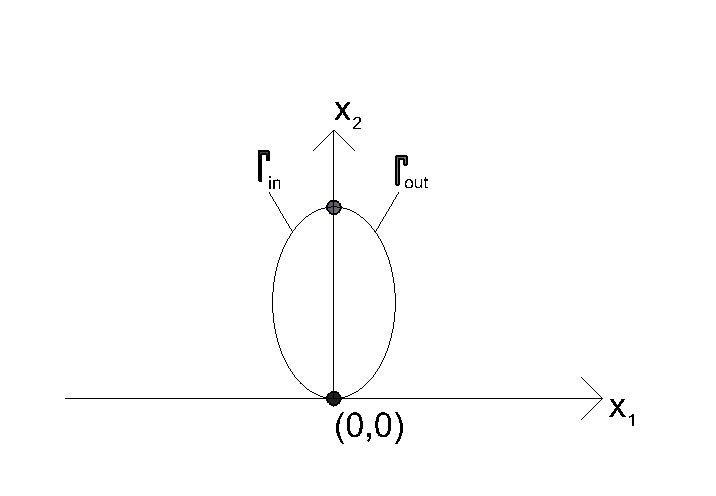}
\caption{The domain. Simple representative case with two singularity points.}
\end{center}
\end{figure}
We assume further that these are the only singularity points, that is, the are no singularity points 
'inside' $\Gamma_{in}$ and $\Gamma_{out}$. 
An example of such domain is shown in Figure 1.
It is well known and was already mentioned in the introduction that existence of regular solutions
to inflow problem \eqref{main_system} requires certain assumptions on the shape of the boundary around the singularity points.
We also need an assumption of this kind, in order to formulate it notice that around each singularity point
the boundary is given as a function $x_2(x_1)$.
We assume it satisfies the condition
\begin{equation} \label{flat}  
\exists N \in \mathbb{N}: \;  |x_2(x_1) - x_2(y_1)| \geq C |x_1-y_1|^N \quad \forall x_1,y_1:\; 
 (x_1,x_2(x_1)),(y_1,x_2(y_1)) \in \Gamma^s_{in}.
\end{equation}
We emphasize that the above condition is required only on the inflow part, therefore it  
can be rewritten as 
\begin{equation} \label{flat_inv}
\exists \delta>0: \; |\underline{x_1}(x_2)-\underline{x_1}(y_2)|
\leq C | x_2-y_2 |^\delta \quad \textrm{with} \quad \delta=\frac{1}{N}.
\end{equation}

\begin{rem} \label{rmk:flat}
Notice that if \eqref{flat} is satisfied by some $N^*$ then it is also satisfied for $N>N*$ for 
small $|x_2-y_2|$. This condition means that 
the boundary around a singularity point must be less flat than some polynomial. Moreover we shall emphasize 
that \eqref{flat} is a necessary but not sufficient condition; we also require sufficient global regularity of the boundary
in order to solve an auxilary elliptic problem. In order to avoid additional technicalities  
we assume that $\Omega$ is a $C^3$ domain. Such regularity is obviously not assured by \eqref{flat},
in particular a Lipschitz boundary may satisfy \eqref{flat} for any $N \geq 1$.  
\end{rem}

\begin{rem} \label{rmk:flat2}
A similar domain is considered in \cite{Kw1} and 
it is worth to compare the condition \eqref{flat_inv} with a similar constraint there with 
$\delta \geq \frac{1}{2}$. Therefore our condition is clearly weaker 
and means that the boundary around 
the singularity points is less flat then some polynomial. 
\end{rem}

\noindent
Another interesting observation concerning the condition \eqref{flat} is that a polynomial 
behaviour of the boundary near the singularity points turns out important in a completely 
different context of singularity of boundary layers in the stationary Munk equation, see \cite{DSR}. 

%After showing the proofs for a representative domain
%described above we explain briefly how they are generalized to wider class of domains.  

Although condition \eqref{flat_inv} seems quite technical in the above formulation,
it is in fact satisfied by a wide class of functions, in particular by piecewise analytical boundaries
what is shown in the following lemma:
\begin{lem}
Assume that $x_2$ is an analytic function of $x_1$ around the singularity points. 
Then (\ref{flat}) holds.
\end{lem}
   
\emph{Proof.} 
By \eqref{def_sing2} we have $x_2'(x_1)=0$ at the singularity points. Therefore 
it is enough to show that if $f: \mathbb{R} \to \mathbb{R}$ is analytic in some $[-r,r]$,
$f(0)=f'(0)=0$ and $f \neq 0$ then
\begin{equation}
|f(x)| \geq C |x|^N \quad \textrm{for} \quad x \in (-l,l)
\end{equation}   
for some $C>0$, $N \geq 2$ and $l<r$ sufficiently small.  
Since $f \neq 0$ and $f$ is analytic, we must have $f^{(n)}(0) \neq 0$ for some $n \geq 2$.
Let $f^{(k)}(0)$ be the first derivative not vanishing in $0$. Then we have
$$
f(x) = \frac{f^{(k)}(0)}{k!} x^k + R^{k+1}(x),
$$
where $|R^{k+1}(x)| \leq M |x|^{k+1}$ for $x \in (-r,r)$.
Hence 
$$
|f(x)| \geq \Big|\frac{f^k(0)}{k!}\Big| |x|^k-M|x|^{k+1} = \Big( \frac{f^k(0)}{k!} - M|x|\Big) |x|^k. 
$$
\begin{flushright} $\square$ \end{flushright}

\subsection{The domain. General case}
Our result holds for a wide class of domains where the inflow and
outflow parts are defined in a natural way. A general setting is to consider inflow and outflow described as
\begin{equation} \label{def_in_out}
\begin{array}{c}
\Gamma_{in}=\{(\underline{x_1}(x_2),x_2): x_2 \in (a,b)\},\\[5pt]
\Gamma_{out}=\{(\overline{x_1}(x_2),x_2): x_2 \in (a,b)\},\\ 
\end{array}
\end{equation}
with singularity points given by 
$\{(\ux(x_2),x_2),x_2=k_{in}^i\}$ and $\{(\ox(x_2),x_2),x_2=k_{out}^j\}$ where %todo figure  
$$
{\rm lim}_{x_2 \to k_{in}^{i-}} |\ux'(x_2)| = {\rm lim}_{x_2 \to k_{in}^{i+} } |\ux'(x_2)| = \infty,
$$$$
{\rm lim}_{x_2 \to k_{out}^{j-}} |\ox'(x_2)| = {\rm lim}_{x_2 \to k_{out}^{j+}}  |\ox'(x_2)|  = \infty.
$$
%In particular it is possible to introduce $\Gamma_0$ of positive measure.  
%Then $\Gamma_{in}$ and $\Gamma_{out}$ are still given by \eqref{def_in_out} 
%but this time we admit jump discontinuities of the function $\ux(x_2)$ and $\ox(x_2)$
%and $\Gamma_0$ is composed of straight intervals connecting these discontinuities (instead 
%of writing a formula for $\Gamma_0$ let us refer to Figure todo ).    
%and \eqref{def_sing2} but with
Furthermore we assume 
\begin{equation} %\label{def_sing3}
\begin{array}{c}
{\rm lim}_{x_2 \to a^+} \underline{x_1}(x_2) = c_1, \quad {\rm lim}_{x_2 \to b^-} \underline{x_1}(x_2) = c_2,\\[5pt]
{\rm lim}_{x_2 \to a^+} \overline{x_1}(x_2) = d_1, \quad {\rm lim}_{x_2 \to b^-} \overline{x_1}(x_2) = d_2
\end{array}
\end{equation}
with $c_i \leq d_i$. Then we have 
\begin{equation}
\Gamma_0 = (c_1,d_1) \times \{a\} \cup (c_2,d_2) \times \{b\}.
\end{equation}
\begin{figure} 
\begin{center}
\includegraphics[width=0.7\textwidth]{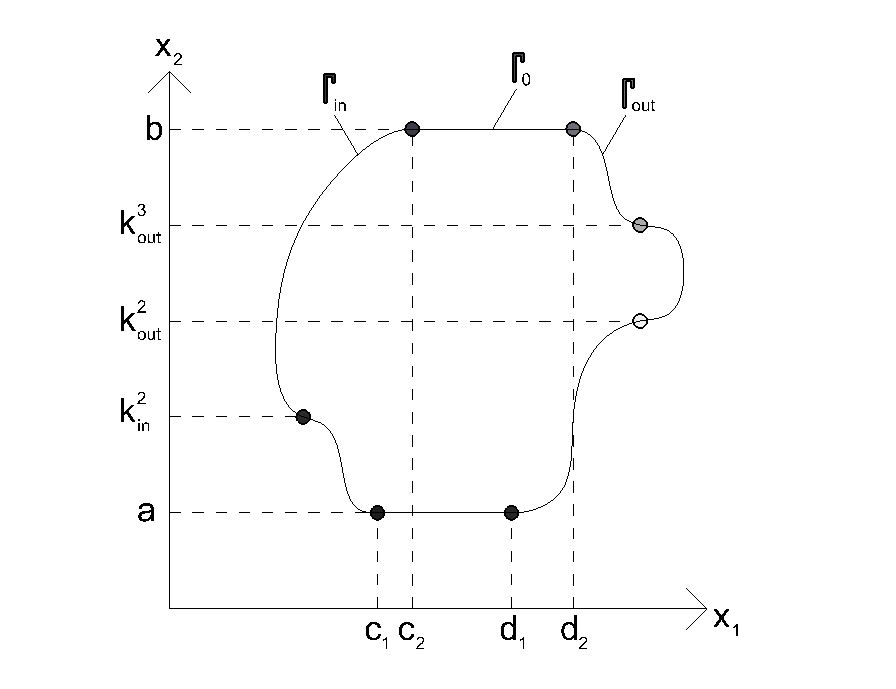}
\caption{The domain. General case}
\end{center}
\end{figure}
An example of above described general domain is shown in Figure 2.
Again, around singularity points $x_2$ is given as a function $x_2^i(x_1)$
and we assume these functions to satisfy \eqref{flat}.

%After showing the proofs for a representative domain
%described above we explain briefly how they are generalized to wider class of domains.  

\subsection{Main result}

We are now in a position to formulate our main result.
\begin{tw} \label{t1} 
Assume that $\Omega$ satifies the conditions introduced above, in particular \eqref{flat_inv}.  
Assume the following regularity of the boundary data:
\begin{align}
&b \in W^{s-1/p}_p(\Gamma), \quad d \in W^{1+s-1/p}_p(\Gamma),\\
&\rho_{in} \in W^s_p(\Gamma_{in}) \cap W^r_p(\Gamma^s_{in}),  \label{c:rhoin}
\end{align}
where 
\begin{equation} \label{c:sp}
s<\delta, \quad r>\frac{s}{\delta}. %\quad p>{\rm max}\{\frac{2}{s}, \frac{\delta}{\delta r-s}\}.
\end{equation}  
Let the pressure be in the form 
\begin{equation} \label{plin}
\pi(\rho)=K\rho
\end{equation}
for some positive constant $K$. Let the viscosity $\mu$ be sufficiently
large compared to $|\Omega|$, $\|\bar v\|_{L_\infty}$ and $K$.
Assume further that the boundary data satisfy following additional assumption: 
\begin{displaymath} %\label{sing}
|(d-n_1)| \sqrt{1+|\underline{x_1}'(x_2)|^2} \leq \kappa <1,  \qquad (A1)  
\end{displaymath}
Assume that $D_0$ defined in (\ref{D0}) is small enough and
let $f$ be large enough on $\Gamma_{in}$. 
%Assume additionally that 
%\begin{displaymath}
%dist(\Gamma_{in},\Gamma_{out}) \geq \lambda >0. \qquad (A3)
%\end{displaymath}
%
Then there exists 
a solution $(v,\rho) \in W^{1+s}_p \times W^s_p$ to the system (\ref{main_system}) such that 
\begin{equation} \label{est_main}
\|v - \bar v\|_{W^{1+s}_p} + \|\rho - \bar \rho\|_{W^s_p} \leq E(D_0),
\end{equation}
where $E(\cdot)$ is a Lipschitz function, E(0)=0.
Moreover, this solution is unique in the class of solutions satisfying (\ref{est_main}).
\end{tw}
\begin{rem} \label{rmk:4}
Condition \eqref{c:rhoin} is required to obtain the estimate for the density near the singularity 
points, the details are shown in the proofs of Lemmas \ref{L6} and \ref{L7}. Therefore, for $\delta$ given by the geometry of the domain in \eqref{flat_inv} 
we can choose arbitrarily small $r$ provided we choose 
sufficiently small $s$, but we have to compensate it with sufficiently large $p$. 
Note that in particular we can assume 
\begin{equation} \label{c:rhoin2}  
\rho_{in} \in W^s_p(\Gamma_{in}) \cap W^1_p(\Gamma^s_{in}). %\quad p>{\rm max}\{ \frac{2}{s},\frac{\delta}{\delta-s} \}.
\end{equation}
\end{rem}
The rest of the paper is organized as follows. In the remaining of the present section we reformulate
the problem \eqref{main_system} introducing perturbations as new unknowns, obtaining system \eqref{system}. 
Then we recall basic properties of the functional spaces
we use. In Section 2 we introduce linearization of \eqref{system} and show a priori estimates.
First we show the energy estimate. Then in Section 2.2 we move to the estimate in $W^s_p$
for the steady transport equation which is the main difficulty in the proof.
This result makes it possible to conclude the estimate in $W^{1+s}_p \times W^s_p$ in Sections 2.3 and 2.4. 
In Section 3 we prove Theorem with an interative scheme using our estimates for the linear system 
to show the convergence of the sequence of approximations. Finally we finish with a short concluding section.     
Without loss of generality we assume in the proofs $v^*=1$ in \eqref{const} except the proof 
of Proposition 3 where we need to track the dependence of the viscosity on this constant.

To remove inhomogeneity from the boundary condition (\ref{system})$_4$, we construct $u_0 \in W^{1+s}_p$ 
such that 
\begin{equation} \label{def_u0}
u_0 \cdot n|_{\Gamma} = d - n_1, \quad
\|u_0\|_{W^{1+s}_p} \leq C \|d - n_1\|_{W^{1+s-1/p}_p(\Gamma)}. 
\end{equation}
For our purpose we find $u_0$ as a solution to the problem
\begin{equation} \label{def_ext}
-\mu \Delta u_0 - (\mu+\nu) \nabla {\rm div} u_0 = 0, \quad 
u_0 \cdot n = d - n_1.
\end{equation}
In order to construct $u_0$ we supply \eqref{def_ext} with condition $u_0 \cdot \tau=0$
and define $u_0 = u_0^1 + u_0^2$, where $u_0^1 \in W^{1+s}_p$ is any extension of the boundary data $d-n_1$
and $u_0^2$ solves
$$
-\mu \Delta u_0^2 - (\mu+\nu) \nabla {\rm div} u_0^2 = \mu \Delta u_0^1 + (\mu+\nu) \nabla {\rm div} u_0^1, \quad
u_0^2|_{\Gamma}=0.
$$
Introducing the perturbations 
$$
u = v - \bar v - u_0 \quad {\rm and} \quad w = \rho - \bar \rho  
$$
we obtain the system
\begin{equation} \label{system}
\begin{array}{lcr}
\partial_{x_1} u -\mu \Delta  u - (\nu + \mu) \nabla {\rm div}\,  u +
K \, \nabla  w =  F(u,w) & \mbox{in} & \Omega,\\[3pt]
{\rm div} \, u + \partial_{x_1}w + (u+u_0) \cdot \nabla  w = G(u,w)
& \mbox{in}& \Omega,\\[3pt]
n\cdot 2\mu {\bf D}( u)\cdot \tau +f \ u \cdot \tau = B
&\mbox{on} & \Gamma, \\[3pt]
n\cdot  u = 0 & \mbox{on} & \Gamma,\\[3pt]
w=w_{in} & \mbox{on} & \Gamma_{in},
\end{array}
\end{equation}
where %$\gamma=\pi'(1)$ 
\begin{align}  \label{FG}
%\begin{array}{c}
F(u,w) &= - w \, (u + \bar v + u_0) \cdot \nabla (u+u_0)-\d_{x_1}u_0 
- (u + u_0) \cdot \nabla (u+u_0),
%+ \mu \Delta u_0 + (\nu+\mu) \nabla \div  u_0,
\nonumber\\[3pt]
G(u,w) &= -(w + 1) \, \div u_0 - w \, \div u,\nonumber\\[3pt]
B &= b - f \tau^{(1)} - 2 \mu \, n \cdot \mathbf{D}(u_0) \cdot \tau.
%\end{array}
\end{align}
Notice that in general $F(u,w)$ contains a term $\mu \Delta u_0 + (\nu+\mu) \nabla {\rm div}\, u_0$
which vanishes due to our definition of $u_0$. 
%
%\begin{equation} \label{f1_exp}
%F_1(w)=\pi(w+1) - \pi'(1)w - \pi(1)= 
%w^2 \int_0^1 \pi''(1-\theta+\theta w)d\theta.
%\end{equation} 
%Notice that we can modify $F_1$ by any constant, therefore we subtract $\pi(1)$ to have \eqref{f1_exp}.
It can be seen easily that 
\begin{align}  \label{a7}
%\begin{array}{c}
\|F(u,w)\|_{L_p} &\leq C [ \|w\|_{W^s_p}\|u\|_{W^{1+s}_p} + E(\|\nabla u,u,w\|_{L_p}+1) ],\nonumber\\[3pt]
%\|F_1(w)\|_{W^s_p} \leq C \|w\|_{W^s_p}^2, \\[5pt]
\|G(u,w)\|_{W^s_p} &\leq C [  \|w\|_{W^s_p}\|u\|_{W^{1+s}_p} + E(\|w\|_{W^s_p}+1) ],\nonumber\\[3pt]  
\|B\|_{W^{s-1/p}_p} &\leq C [\|b-f \tau\|_{W^{s-1/p}_p} + \|d-n_1\|_{W^{1+s-1/p}_p(\Gamma)}],
%\end{array}
\end{align}
where $E$ is a small constant dependent on $\|u_0\|_{W^{1+s}_p}$ and
$V^*$ is a dual space to $V$ defined in \eqref{def_V}.
%The bound on $B$ follows from the trace theorem in $W^s_p$. 

From now on we focus on the system (\ref{system}). Our goal is to show the existence of a solution
$(u,w) \in W^{1+s}_p \times W^s_p$ for given small function $u_0 \in W^{1+s}_p$.
Recall in particular that solution to the system \eqref{system} is used in the definition of solution
to the the original problem \eqref{main_system}. 
In order to define the solution to \eqref{system} we need its weak formulation. 
For this purpose we apply the identity 
\begin{align} \label{basic_id}
%\begin{array}{c}
&\int_{\Omega} (-\mu \Delta u - (\nu +\mu)\nabla \, {\rm div} u) \cdot v \, dx =
\int_{\Omega} 2 \mu {\bf D}(u): \nabla \,v + \nu \, {\rm div} \,u \, {\rm div} \,v \, dx \nonumber\\[3pt]
&- \int_{\Gamma} n \cdot [2 \mu {\bf D}(u)] \cdot v \, d\sigma - \int_{\Gamma} n \cdot [\nu ({\rm div} u) {\bf{Id}} ] \cdot v \, d\sigma.
%\end{array}
\end{align}
Then a natural definition is following
\begin{df} \label{def_sol1}
A regular $W^s_p$ solution to the problem \eqref{system} is a couple $(u,w) \in W^{1+s}_p \times W^s_p$
such that
\begin{align} \label{weak1}
%\begin{array}{c} 
&\int_{\Omega} \psi \d_{x_1}u \,dx + \int_{\Omega} (2\mu {\bf D}(u):\nabla \psi + \nu \div u \div \psi )\,dx
+\int_{\Gamma} [f(u \cdot \tau)(\psi \cdot \tau)-b(\psi \cdot \tau)]\,d\sigma \nonumber\\
&- \int_{\Omega} w \div \psi \,dx
=\int_{\Omega} F(u,w) \cdot \psi \,dx
\quad \forall \psi \in V 
%\end{array}
\end{align}
and
\begin{align} \label{weak2} 
%\begin{array}{c}
&\int_{\Omega}\phi \div u \,dx - \int_{\Omega} w ((\bar v+u+u_0) \cdot \nabla \phi + \phi \div (u+u_0))\,dx
+ \int_{\Gamma_{in}} w_{in}\phi \,d\sigma = \nonumber\\
&= \int_{\Omega} G(u,w)\phi\,dx 
\quad \forall \phi \in C^1(\Omega): \phi|_{\Gamma \setminus \Gamma_{in}}=0.
%\end{array}
\end{align}

\end{df}

\section{Linearization and a priori bounds}
In this section we derive \emph{a priori} estimates for the following
linearization of system \eqref{system}
\begin{equation} \label{system_lin}
\begin{array}{lcr}
\partial_{x_1} u -\mu \Delta  u - (\nu + \mu) \nabla {\rm div}\,  u +
K \, \nabla  w = F & \mbox{in} & \Omega,\\[3pt]
{\rm div} \, u + \partial_{x_1}w + U \cdot \nabla  w = G
& \mbox{in}& \Omega,\\[3pt]
n\cdot 2\mu {\bf D}( u)\cdot \tau +f \ u \cdot \tau = B
&\mbox{on} & \Gamma, \\[3pt]
n\cdot  u = 0 & \mbox{on} & \Gamma,\\[3pt]
w=w_{in} & \mbox{on} & \Gamma_{in},\\[3pt]
\end{array}
\end{equation}
where $U \in W^{1+s}_p$ is small and satisfies $U \cdot n|_{\Gamma}=d-n_1$
and on the rhs we have 
\begin{equation} \label{reg_rhs}
F \in L_p(\Omega), \; 
G \in W^s_p(\Omega) \; \textrm{and} \; B \in W^{s-1/p}_p(\Gamma).
\end{equation}
We start with the energy estimate, then we deal with the steady transport equation
which is the main difficulty of our proof and finally we show the estimate in the solution
space.
\subsection{Energy estimate}
In this section we show energy estimate for the solutions of \eqref{system_lin}.
\begin{lem}
Let $(u,w)$ be a sufficiently smooth solution to system (\ref{system_lin}) 
with given functions $(F,G,B) \in (V^* \times L_2 \times L_2(\Gamma))$. Then 
\begin{equation} \label{ene}
\|u\|_{W^1_2}+\|w\|_{L_2} \leq C \big[ \|F\|_{V^*} + \|G\|_{L_2} + \|B\|_{L^2(\partial \Omega)} \big],
\end{equation}
where $C$ is independent from the boundary data. 
\end{lem}
%
%In order to show \eqref{ene} we apply a standard energy method.
\noindent
\emph{Proof.}
Multiplying the first equation of (\ref{system_lin}) by $u$ and integrating over $\Omega$ we
get using the boundary condition \eqref{system_lin}$_3$ and the identity \eqref{basic_id}:
\begin{align} \label{a5}
& \int_\Omega 2\mu{\bf D}^2(u) +\nu \div^2udx + \int_{\partial \Omega} \left(f + \frac{n_1}{2}\right) u^2 d \sigma 
+ \int_\Omega K \nabla w u dx = \nonumber \\
& -\int_\Omega u \, \partial_{x_1}
u \,  dx +\int_{\Omega} F u dx + \int_{\partial \Omega} B(u \cdot \tau) \, d\sigma.
\end{align}

The boundary term on the lhs will be positive for $f\geq0$ on $\Gamma_{out}$ and $f \geq \frac{n_1}{2}$
on $\Gamma_{in}$.
Next we integrate by parts the last term of the l.h.s of (\ref{a5}).
Using (\ref{system_lin})$_2$ we obtain:
\begin{align} \label{a6}
K \int_{\Omega} u\nabla w  dx &= - K \int_{\Omega} w\div u  dx=
K \int_{\Omega} (w\partial_{x_1} w + U \cdot w\nabla w  - G w )dx \nonumber\\
&=K \Big[ \frac{1}{2}\int_{\Gamma} w^2 n_1 d\sigma - \frac{1}{2}\int_{\Omega} w^2{\rm div}\,U - \int G w dx \Big].
\end{align}
We will also use the following Korn inequality:
\begin{equation} \label{korn}
\int_{\Omega} 2 \mu {\bf D^2}(u) + \nu \div^2 u \,dx \geq C_K \| u \|_{W_2^1(\Omega)}^2,
\end{equation}
where $C_K=C_K(\Omega)$.
Using \eqref{a5}, \eqref{a6} and \eqref{korn} we get:
\begin{align*}
&\| u \|_{W_2^1(\Omega)}^2 + \int_{\Gamma_{out}} w^2 n_1 \, d\sigma  \leq\\
&\leq C \left[ \int_{\Omega} w^2 {\rm div} \, U \, dx + \gamma \int_{\Omega} Gw \,dx + \int_{\Omega} Fu \,dx
+ \int_{\Gamma} B (u \cdot \tau) \, d\sigma + \int_{\Gamma_{in}} w_{in}^2 n_1 d\sigma \right].
\end{align*}
Next, using H\"older and Young inequalities,
the fact that $w^2 n_1>0$ on $\Gamma_{out}$
and the trace theorem to the boundary term we get for any $\delta>0$:
\begin{displaymath} 
\begin{array}{c}
\|u\|^2_{W^1_2(\Omega)} \leq \big( \delta + \|U\|_{W^{1+s}_p(\Omega)}\big) \|w\|_{L_2(\Omega)}^2 + C(\delta) \|G\|_{L_2(\Omega)}^2
+C \left[ \big( \|F\|_{V^*} + \|B\|_{L_2(\Gamma)} \big) \|u\|_{W^1_2(\Omega)} \right],
\end{array}
\end{displaymath}
which yields
\begin{equation} \label{u1} 
\begin{array}{c}
\|u\|_{W^1_2(\Omega)} \leq \big( \delta + \|U\|_{W^{1+s}_p(\Omega)}\big) \|w\|_{L_2(\Omega)} +C(\delta)\|G\|_{L_2(\Omega)}+C( \|F\|_{V^*}+\|B\|_{L_2(\Gamma)} ).
\end{array}
\end{equation}
To estimate the first term of the r.h.s. we find a
bound on $\|w\|_{L_2}$.
%Let us define:
%$$
%a = min \, \{ \underline{x_1}(x_2): \, x_2 \in (0,b) \}
%$$
From $(\ref{system_lin})_2$ we have
$$
\partial_{x_1}w + U\cdot \nabla w = G - {\rm div} \,u.
$$
In order to estimate $\|w\|_{L_2}$, for $x \in \Omega$ let us denote by $\gamma_x$ a characteristic
of the operator $\partial_{x_1}+U \cdot \nabla$ connecting $x$ with $\Gamma_{in}$, i.e. solution to
\begin{equation} \label{def:gamma}
\left\{ \begin{array}{c}
\dot{\gamma_x}(t)=-(1-U_1(\gamma_x(t)),U_2(\gamma_x(t))),\\ 
\gamma_x(0)=x.
\end{array} \right.
\end{equation}
Notice that we take the tangential vector with opposite sign since we consider a backwards 
characteristic starting at a given point inside $\Omega$.
Denote by $x^{in}(x)$ the intersection of $\gamma_x$ with $\Gamma_{in}$.
Due to regularity and smallness of $U$, $\gamma_x$ is close to a straight line $\{x_2=c\}$. Now we can write
\begin{equation} \label{wx}
w(x) = w_{in}(x^{in}(x)) + \int_{\gamma_x} (G - \rm{div}\,u) \,dl_{\gamma_x}.
\end{equation}
By Jensen inequality we have
$$
\left( \int_\gamma (G - \rm{div}\,u) dl_{\gamma_x} \right) ^2 \leq |\gamma_x| \int_{\gamma_x} |G|^2+|{\rm div}\,u|^2 dl_{\gamma_x}.
$$
Hence applying \eqref{wx} we get
\begin{equation} \label{wl2}
\|w\|_{L_2(\Omega)} \leq C(\Omega) \left( \|w_{in}\|_{L_2(\Gamma_{in})} + \|G\|_{L_2(\Omega)} + \|u\|_{W^1_2(\Omega)} \right). 
\end{equation}
Now we combine \eqref{u1} and \eqref{wl2}. By smallness of $U$ we can fix $\delta$ in \eqref{u1}
small enough to put the term $\delta+\|U\|_{W^{1+s}_p}$ on the left obtaining \eqref{ene},
which completes the proof of the lemma. 
\qed

\subsection{Steady transport equation}
In this section we show the estimate in $W^s_p$ for the steady transport equation with inflow condition
which is a crucial step in showing a priori estimate for the linear problem \eqref{system_lin}.
\begin{prop} \label{lem_trans}
Let $\Omega$ be a set defined at the end of Section 1, satisfying \eqref{flat_inv}.
Let $w$ solve
\begin{equation} \label{trans}
\bar K w+w_{x_1}+ U \cdot \nabla w = H, \quad w|_{\Gamma_{in}}=w_{in},
\end{equation}
where $\bar K$ is a positive constant,  $U \in W^1_\infty(\Omega)$ is small and satisfies $U \cdot n|_{\Gamma}=d-n_1$ satisfying (A1).
Assume $H \in W^s_p(\Omega)$ and let $w_{in}$, $s$ and $p$ 
satisfy the assumptions of Theorem \ref{t1}.
Then 
\begin{equation} \label{est_trans}
\|w\|_{W^s_p(\Omega)} \leq C [\|H\|_{W^s_p(\Omega)} + \|w_{in}\|_{W^s_p(\Gamma_{in}) \cap W^r_p(\Gamma_s)}], 
\end{equation} 
where $C=C(s,p,\Omega)$ and $r$ is from \eqref{c:sp}.
\end{prop}
\begin{rem}
Notice that the assumption $U \in W^1_\infty(\Omega)$ is weaker than $U \in W^{1+s}_p(\Omega)$ for $sp>2$
due to the imbedding theorem.   
\end{rem}
\begin{rem}
Proposition \ref{lem_trans} generalizes the result from \cite{P3} where it is assumed $U=[u^1,0]$ which enables 
much simpler proof that the one presented below and no additional regularity of the
density around the singularity points is required.     
\end{rem}
\noindent
For simplicity we set $\bar K=1$ which is allowed as we assume anyway sufficient smallness of the data.
The proof of Proposition \ref{lem_trans} is quite technical since we have to treat carefully the 
boundary terms. For the reader's convenience we divide it into several lemmas.  
\begin{lem}
Let the assuptions of Proposition \ref{lem_trans} hold. Then 
\begin{multline} \label{est_trans_1}
\|w\|_{W^s_p(\Omega)} + \frac{1}{p}\iint_{\Omega^2} \frac{\partial_{x_1}|w(x)-w(y)|^p+\partial_{y_1}|w(x)-w(y)|^p}{\phi_{\epsilon}(x,y)}dxdy \\
+ \frac{1}{p}\iint_{\Omega^2} \frac{ U(x)\cdot \nabla_x|w(x)-w(y)|^p + U(y)\cdot \nabla_y|w(x)-w(y)|^p }{\phi_{\epsilon}(x,y)}dxdy 
\leq \|H\|_{W^s_p(\Omega)} \|w\|_{W^s_p}^{p-1}.
\end{multline}
\end{lem}
\emph{Proof.} 
Recalling the definition of Sobolev-Slobodetskii norm we write \eqref{trans} in $x$ and $y$.
Using identities of a kind of $\nabla_x w(y) = 0$ we can write
$$
w(x)+ \partial_{x_1}[w(x)-w(y)] +U(x) \cdot \nabla_x [w(x) - w(y)] = H(x), 
$$$$
w(y)+ \partial_{y_1} [w(y)-w(x)] +U(y) \cdot \nabla_y [w(y) - w(x)] = H(y).
$$
We multiply the first equation by $\frac{|w(x)-w(y)|^{p-2}(w(x)-w(y))}{\phi_{\epsilon}(x,y)}$
and the second by \newline
$\frac{|w(x)-w(y)|^{p-2}(w(y)-w(x))}{\phi_{\epsilon}(x,y)}$,
where
\begin{equation*}
\phi_{\epsilon}(x,y)=\epsilon+|x-y|^{2+sp}. 
\end{equation*}
Then we add the equations
and integrate twice over $\Omega$, w.r.t. $x$ and $y$.
Since 
$$
w(x)(w(x)-w(y))+w(y)(w(y)-w(x))=[w(x)-w(y)]^2
$$
and
$$
H(x)(w(x)-w(y))+H(y)(w(y)-w(x))=(H(x)-H(y))(w(x)-w(y)),
$$
we obtain on the left hand side
\begin{equation} \label{e1}
\iint_{\Omega^2} \frac{|w(x)-w(y)|^p}{\phi_{\epsilon}(x,y)} \,dxdy \to_{\epsilon \to 0} \|w\|_{W^s_p}^p.
\end{equation}
On the r.h.s we have using H\"older inequality:
\begin{align} \label{e2}
&\iint_{\Omega^2} \frac{(H(x)-H(y))|w(x)-w(y)|^{p-2}(w(x)-w(y))}{\phi_{\epsilon}(x,y)} \leq
\nonumber\\
&\Big( \iint_{\Omega^2}  \frac{ |H(x)-H(y)|^p }{\phi_{\epsilon}(x,y)} \Big)^{1/p}
\Big( \iint_{\Omega^2}  \frac{ |w(x)-w(y)|^p }{\phi_{\epsilon}(x,y)} \Big)^{1-1/p} 
\to_{\epsilon \to 0} \|H\|_{W^s_p} \|w\|_{W^s_p}^{p-1},
\end{align}
which completes the proof of \eqref{est_trans_1}. \qed

The rest of the proof of Proposition \ref{lem_trans}
consist in dealing with the integral terms in \eqref{est_trans_1}. This is where
all the difficulties are hidden and our assumptions on the boundary and boundary data will come into play.
We start with observing that under the assumptions of Proposition \ref{lem_trans} we have 
\begin{equation} \label{out:pos}
n_1+U\cdot n \geq 0 \quad {\rm on} \quad \Gamma_{out}.
\end{equation}
%\emph{Proof.} Let us define $\Gamma_{out}=\Gamma_{out}^s \cup \Gamma_{out}^r$ 
%where $\Gamma_{out}^s$ is lying close enough to singularity points that 
%\begin{equation} \label{def:out_s}
%\sqrt{1+[\ox'(x_2)]^2}\leq (1+\frac{1-\kappa}{2\kappa}|\ox'(x_2)| \quad {\rm for} \quad 
%x_2,\ox(x_2) \in \Gamma_{out}^s. 
%\end{equation}
Indeed, we have 
\begin{equation} \label{normal}
n(x)|_{\Gamma_{out}}=\frac{[1,-\overline{x_1}'(x_2)]}{\sqrt{1+[\overline{x_1}'(x_2)]^2}}, \qquad
n(x)|_{\Gamma_{in}}=\frac{[-1,\underline{x_1}'(x_2)]}{\sqrt{1+[\underline{x_1}'(x_2)]^2}},
\end{equation}
therefore by (A1)  
$$
\left|\frac{U \cdot n}{n_1}\right|=\sqrt{1+[\ox'(x_2)]^2}|U\cdot n|\leq \kappa, 
$$
and so \eqref{out:pos} holds. 
%Next we show a bound from below on $n_1$ on $\Gamma_{out}^r$. For this purpose denote 
%$\sigma(\kappa)=\frac{1-\kappa}{2\kappa}$. Then \eqref{def:out_s} implies
%\begin{equation}
%|\ox'(x_2)| \geq \frac{1}{\sqrt{\sigma(\kappa)(2+\sigma(\kappa))}} \quad {\rm on} \quad \Gamma_{out}^s
%\end{equation}
%and so by \eqref{normal} 
%\begin{equation}
%n_1 > \frac{1}{1+\sigma(\kappa)|\ox'(x_2)|} \geq \frac{1}{(1+\sigma(\kappa))\sqrt{\sigma(\kappa)(2+\sigma(\kappa))}}=:A(\kappa).
%\end{equation}
%Notice that $A(\kappa)$ is increasing in $\kappa$. Moreover, if (A1) holds with some $\kappa<\frac{1}{2}$
%then it holds with $\kappa = \frac{1}{2}$ and if $\kappa > \frac{1}{2}$ then $A(\kappa)>A(\frac{1}{2})$
%Therefore \eqref{out:pos} holds on $\Gamma_{out}^r$ provided $|U \cdot n| < A(1/2)$.
%\qed

We now proceed with transforming \eqref{est_trans_1}. 
Since we want to have a $W^s_p$ norm, we have to get rid of the derivatives of $w$ and for this purpose we 
integrate by parts. This step is presented in the following lemma  
\begin{lem}
Let the assumptions of Proposition \ref{lem_trans} be satisfied. Then
\begin{align} \label{est_trans_2}
&\|w\|_{W^s_p(\Omega)}^p + \frac{2}{p} \int_{\Omega} \int_{\Gamma_{in}} \frac{|w(x)-w(y)|^p (n_1+U \cdot n)(x)}{\phi_{\epsilon}} dS(x)dy \leq \nonumber\\[3pt]
& \leq \|H\|_{W^s_p(\Omega)} \|w\|_{W^s_p(\Omega)}^{p-1} + C \|\nabla U\|_{L_{\infty}(\Omega)} \|w\|_{W^s_p(\Omega)}^p.
\end{align}
\end{lem}
\emph{Proof.} We will obtain \eqref{est_trans_2} from \eqref{est_trans_1}.
Let us start with the first integral term on the lhs of \eqref{est_trans_1}.
We have
\begin{equation} \label{wx1_1}
 \frac{\partial_{x_1}|w(x)-w(y)|^p}{\phi_{\epsilon}(x,y)}=
\partial_{x_1} \Big( \frac{|w(x)-w(y)|^p}{\phi_{\epsilon}(x,y)} \Big) - |w(x)-w(y)|^p \partial_{x_1} 
\Big( \frac{1}{\phi_{\epsilon}(x,y)} \Big).
\end{equation}
By the definition of $\phi_{\epsilon}(x,y)$ we have
%$$
%\partial_{x_i} \frac{1}{\phi_{\epsilon}(x,y)} = - \frac{ (2+sp) \, |x-y|^{sp} (x_i-y_i) }{\phi_{\epsilon}^2}
%= - \partial_{y_i} \frac{1}{\phi_{\epsilon}(x,y)}.
%$$
%In particular,
%
\begin{equation} \label{phi_sym}
\nabla_x\phi_{\epsilon}(x,y)=-\nabla_y\phi_{\epsilon}(x,y),
\end{equation}
therefore using \eqref{wx1_1} we get
\begin{align} \label{wx1_2} 
&\iint_{\Omega^2} \Big\{ \frac{\partial_{x_1}|w(x)-w(y)|^p}{\phi_{\epsilon}(x,y)}
+ \frac{\partial_{y_1}|w(y)-w(x)|^p}{\phi_{\epsilon}(x,y)} \Big\} \,dxdy = 
%$$$$
%= \frac{1}{p} \iint_{\Omega^2} \frac{\partial_{x_1}|w(x)-w(y)|^p}{\phi_{\epsilon}(x,y)} \,dxdy
%+ \frac{1}{p} \iint_{\Omega^2} \frac{\partial_{y_1}|w(x)-w(y)|^p}{\phi_{\epsilon}(x,y)} \,dxdy =
\nonumber\\[3pt]
&= \iint_{\Omega^2} \Big\{ \partial_{x_1} \Big( \frac{|w(x)-w(y)|^p}{\phi_{\epsilon}(x,y)} \Big)
+  \partial_{y_1} \Big( \frac{|w(x)-w(y)|^p}{\phi_{\epsilon}(x,y)} \Big) \Big\} \,dxdy
\nonumber\\[3pt]
&- \iint_{\Omega^2} \Big\{ |w(x)-w(y)|^p \partial_{x_1} \Big( \frac{1}{\phi_{\epsilon}(x,y)} \Big) 
+ |w(x)-w(y)|^p \partial_{y_1} \Big( \frac{1}{\phi_{\epsilon}(x,y)} \Big) \Big\} \,dxdy. 
\end{align}
Taking into account \eqref{phi_sym}, the integrand in the second integral on the rhs of \eqref{wx1_2} vanishes identically. 
Combining \eqref{phi_sym} with the identities
\begin{equation} \label{nabla_x_sym}
\nabla_x |w(x)-w(y)|^p = p |w(x)-w(y)|^{p-2} (w(x)-w(y)) \nabla_x w(x)
\end{equation}
and
\begin{equation} \label{nabla_y_sym}
\nabla_y |w(x)-w(y)|^p = - p |w(x)-w(y)|^{p-2} (w(x)-w(y)) \nabla_y w(y),
\end{equation}
we see that the first integral on the rhs of \eqref{wx1_2} adds up to
\begin{equation} \label{bd1}
\frac{2}{p} \iint_{\Omega^2} \partial_{x_1} \Big( \frac{|w(x)-w(y)|^p}{\phi_{\epsilon}(x,y)} \Big) \,dxdy = 
\frac{2}{p} \int_{\Omega} \int_{\Gamma} \frac{|w(x)-w(y)|^p n_1}{\phi_{\epsilon}(x,y)} \,dS(x) \,dy.
\end{equation}
Now consider the second integral on the lhs of \eqref{est_trans_1}.
We have
\begin{align} \label{ux}
&\iint_{\Omega^2}U(x) \frac{\nabla_x|w(x)-w(y)|^p}{\phi_{\epsilon}(x,y)} \,dxdy=
\nonumber\\[3pt]
&=\iint_{\Omega^2}U(x) \nabla_x \Big( \frac{ |w(x)-w(y)|^p }{\phi_{\epsilon}(x,y)} \Big)\,dxdy
- \iint_{\Omega^2} U(x) |w(x)-w(y)|^p \nabla_x \Big( \frac{1}{\phi_{\epsilon}(x,y)} \Big)\,dxdy
\end{align}
and
\begin{align} \label{uy} 
&\iint_{\Omega^2}U(y) \frac{\nabla_y|w(x)-w(y)|^p}{\phi_{\epsilon}(x,y)}\,dxdy=
\nonumber\\[3pt]
&=\iint_{\Omega^2}U(y) \nabla_y \Big( \frac{ |w(x)-w(y)|^p }{\phi_{\epsilon}(x,y)} \Big)\,dxdy
-\iint_{\Omega^2} U(y) |w(x)-w(y)|^p \nabla_y \Big( \frac{1}{\phi_{\epsilon}(x,y)} \Big)\,dxdy.
\end{align}
Recalling \eqref{phi_sym} we see that the terms with $\nabla \big( \frac{1}{\phi} \big)$ cancel. But this
time the integrand does not vanish identically and we have to check whether the sum of these integrals
makes sense when $\epsilon \to 0$. This sum equals
$$
\left| \iint_{\Omega^2} [U(x)-U(y)] |w(x)-w(y)|^p \frac{(2+sp) |x-y|^{sp}(x-y) }{ (\epsilon + |x-y|^{2+sp})^2 }\,dxdy \right|
\longrightarrow_{\epsilon \to 0}
$$$$
(2+sp) \left| \iint_{\Omega^2} \frac{U(x)-U(y)}{x-y} \frac{|w(x)-w(y)|^p}{|x-y|^{2+sp}}\,dxdy \right| \leq
C \|\nabla U\|_{L_\infty} \|w\|_{W^s_p}^p.
$$
Now consider the terms with $\nabla \frac{|w(x)-w(y)|^p}{\phi_{\epsilon}}$ on the r.h.s. of \eqref{ux} and \eqref{uy}. 
Combining \eqref{nabla_x_sym} and \eqref{nabla_y_sym} with \eqref{phi_sym} 
we see that these terms add up to
\begin{align} \label{e4_1}
&\frac{2}{p} \iint_{\Omega^2} U(x) \nabla_x \Big( \frac{ |w(x)-w(y)|^p }{\phi_{\epsilon}(x,y)} \Big) \,dxdy=
\nonumber\\[3pt]
&= - \int_{\Omega} dy \int_{\Omega} \frac{|w(x)-w(y)|^p}{\phi_{\epsilon}(x,y)} {\rm div}_x U(x) \,dx
+ \int_{\Omega} dy \int_{\Gamma} \frac{|w(x)-w(y)|^p}{\phi_{\epsilon}(x,y)} U(x) \cdot n \,dS(x).
\end{align}
The first integral is straighforward:
\begin{equation} \label{e4_2}
\Big| \int_{\Omega} dy \int_{\Omega} \frac{|w(x)-w(y)|^p}{\phi_{\epsilon}(x,y)} {\rm div}_x U(x) \,dx \Big| \leq C(\Omega) \|\nabla U\|_{L_{\infty}(\Omega)} \|w\|_{W^s_p(\Omega)}^p.
\end{equation}
Combining \eqref{est_trans_1}, \eqref{bd1}, \eqref{e4_1} and \eqref{e4_2} we get \eqref{est_trans_2}.
%Notice that it is enough to consider these parts of the boundary in the inner integral in \eqref{est_trans_2} since on the rest of the boundary we have
%$n^{(1)}+U \cdot n>0$ by smallness of $U$.

%{\bf Piotr: to sie zmieni jak sie wyjasni to szacowanie tez notacja $n_1$ by sie przydala...}

\qed

In order to complete the proof of Proposition \ref{lem_trans}
it remains to find an estimate on the integral term in \eqref{est_trans_2}.
Notice that it is sufficient to consider the integral over $\Gamma_{in}$ since the outflow 
part is nonnegative due to \eqref{out:pos} and therefore we already skipped it in \eqref{est_trans_2}.
%For this purpose let us first focus on 
%\begin{equation} \label{w0}
%\int_{\Omega} \int_{\Gamma_{in}} \frac{|w(x)-w(y)|^p n_1}{\phi_{\epsilon}(x,y)} dS(x)dy.
%\end{equation}
However, the inflow term is negative because of $n_1$ and therefore must be estimated.
The treatment of this term is quite technical and in fact constitutes the core of the proof of Proposition \ref{lem_trans}. 
Let us start with introducing some further notation.
In the estimate which will follow, $y=(y_1,y_2)$ will denote a point inside $\Omega$ while 
$\underline{x}=(\uxx_1,\uxx_2)$ a point on $\Gamma_{in}$. In several places we will replace an integral with respect to 
the boundary measure on $\Gamma_{in}$ by an integral w.r.t $x_2$. Then we will write 
$\uxx(x_2)=(\ux(x_2),x_2)$. At this point we should also fix $\eta$ in the definition \eqref{def_gamma_s}. 
For this purpose observe that (A1) implies 
\begin{equation} \label{infl:small}
|(d-n_1)| \sqrt{1+|\underline{x_1}'(x_2)|^2} \leq  \frac{\tilde \kappa}{\sqrt{1+\lambda^2}}
\end{equation}
for some $\kappa < \tilde \kappa <1$ and $\lambda>0$. Now in \eqref{def_gamma_s} we choose $\eta$ such that 
\begin{equation} \label{def:lambda}
|x_2'(x_1)|<\lambda \quad {\rm on} \; \Gamma^s_{in}.
\end{equation}
Notice that such $\eta$ exists due to assumed regularity of the boundary. Next we define 
\begin{align} \label{def:omega_in}
&\Gamma_{in}^r:=\{ x \in \Gamma_{in}: \; \eta < x_2 < b-\eta\}, \nonumber\\
&\Omega_{in}^r:=\{ x \in \Omega: \; \eta < x_2 < b-\eta, \quad x_1-\ux(x_2)<\eta \}, \nonumber\\  
&\Omega_{in}^s:=\{ x \in \Omega: \; x_2<2\eta \; \vee \; b-2\eta<x<b\}, \nonumber\\
&\Omega_{in}=\Omega_{in}^r \cup \Omega_{in}^s,
\end{align}
Therefore, $\Gamma_{in}^r$ and $\Omega_{in}^r$ are, respectively, a part of $\Gamma_{in}$ 
and its neighbourhood which are at some fixed distance from singularity points. In particular,
we have 
\begin{equation} \label{reg}
|\ux'(x_2)|\leq C(\eta) \quad {\rm for} \quad (\ux(x_2),x_2) \in \Gamma_{in}^r.
\end{equation}  
Notice also that $\Gamma^s_{in}=\Omega^s_{in} \cap \Gamma_{in}$.
Next, for $y=(y_1,y_2)$ let us denote by $\gamma_y$  
a characteristic curve of the operator $(\partial_{x_1}+U \cdot \nabla)$ connecting $y$ 
with $\Gamma_{in}$, i.e. solution to \eqref{def:gamma} with $\gamma_y(0)=y$.
Furthermore, we denote the intersection of $\gamma_y$ with $\Gamma_{in}$ by $y^{in}(y)$
(see Figure 3).  
%By smallness of $\|U\|_{L_\infty}$, $\gamma_y$ is close to straight line and in particular
%\begin{equation} \label{lgamma}
%|\gamma_y| \simeq |y_1 - \underline{x_1}(y_2)| \simeq |\underline{y}-y|. 
%\end{equation}
%
We will need the following auxiliary result
\begin{lem} \label{L5}
Assume that 
\begin{equation} \label{sing:0}
y \in \Omega_{in}^r, \quad \underline{x} \in \Gamma_{in}^r. 
\end{equation}
Then 
\begin{align} 
|\uxx-\uyy(y)|+|\uyy(y)-y| \leq C|\uxx-y|, \label{sing:1} \\
|\gamma_y| \leq (1+E)|y_1-\ux(y_2)|. \label{sing:1b}
\end{align}
for some constant $C>0$ and a small constant $E$.
\end{lem} 
\emph{Proof.} W.l.o.g. we can restrict to assuming that $x_2$ is a decreasing function of $x_1$ on $\Gamma_{in}$ as the opposite case is analogous.
Then it is convenient to distiguish $3$ cases:
\begin{align} \label{cases}
&\textrm{{\bf (a)}} \quad \uxx_1 \leq \uyy_1(y) \; \textrm{(see Figure 4)}; \nonumber\\
&\textrm{{\bf (b)}} \quad \uyy_1(y) < \uxx_1 < y_1 \; \textrm{(see Figure 3)}; \\
&\textrm{{\bf (c)}} \quad \uxx_1 \geq y_1; \nonumber
\end{align}
In the case {\bf (a)} \eqref{sing:1} is obvious with $C=2$ as $|\uxx-\uyy(y)| \leq |\uxx-y|$ 
and $|\uyy(y)-y| \leq |\uxx-y|$.
The case {\bf (c)} is treated similarly to {\bf (b)}, therefore we don't show it on a separate figure and focus 
on the proof for {\bf (b)}. In that case it is convenient to define 
by $\tilde y(x_2,y)$ the intersection of a line  
$\{ (\underline{x_1}(x_2),t): t \in \mathbb{R} \}$ with a characteristic $\gamma_y$ connecting $\Gamma_{in}$ with $y$.
Such intersection is well defined in the case {\bf b)}, see Figure 3. 
We have for some $\vt \in [\uxx_2,y_2]$
\begin{align} \label{sing:2}
|\uxx-\uyy(y)|\leq & C(|\uxx_1-\uyy_1(y)|+|\uxx_2-\uyy_2(y)|)=(C\ux'(\vt)+1)|\uxx_2-\uyy_2(y)| \nonumber\\ 
\leq & C|\uxx_2-\tilde y_2(x_2,y)| \leq C|\uxx-\tilde y(x_2,y)|\leq C|\uxx-y|,
\end{align}
where we have used \eqref{reg}. Therefore also
\begin{equation*}
|\uyy(y)-y|\leq|\uyy(y)-\uxx|+|\uxx-y|\leq C|\uxx-y|,
\end{equation*} 
which completes the proof of the first assertion. To show the second one 
we can assume w.l.o.g. that $\uyy_2(y) \geq y_2$. We have 
$$
|\uyy_2(y)-y_2|=\int_{\uyy_1(y)}^{y_1} U_2(t,\gamma(t))dt \leq \|U\|_{\infty}|\uyy_1(y)-y_1|,
$$
therefore using again \eqref{reg} we get
$$
|y_1-\uyy_1(y)|\leq |y_1-\ux(y_2)|+|\ux(y_2)-\uyy_1(y)|+|\uyy_2(y)-y_2|\leq 
|y_1-\ux(y_2)|+(|\ux'(\theta)|+1)|\uyy_2(y)-y_2|
$$$$
\leq |y_1-\ux(y_2)|+ (|\ux'(\theta)|+1)\|U_2\|_{\infty}|y_1-\uyy_1(y)|,
$$ 
so for sufficiently small $\|U\|_{\infty}$ we obtain 
$$
|y_1-\uyy_1(y)|\leq (1+E)|y_1-\ux(y_2)|,
$$
and to complete the proof it is enough to note that
$$
|\gamma_y|=\int_{\uyy_1(y)}^{y_1}\sqrt{[1+U_1(\gamma_y(s))]^2+(U_2(\gamma_y(s)))^2} ds \leq (1+E) |y_1-\uyy_1(y)|.
$$

\qed

We shall emphasize that \eqref{sing:1} holds except a neighborhood of singularity points. 
The point is that $|\ux'(x_2)|$ is unbounded as we approach the singularity. Therefore, as 
\begin{equation} \label{sing:18}
|\uxx(x_2)-\uyy(y_2)|\geq |\ux(x_2)-y^{in}_1(y_2)|=|\ux'(\vt)||x_2-y_2|
\end{equation}
for some $\vt \in [x_2 , y_2]$, we no longer have \eqref{sing:1}.
In order to control $|\uxx-y^{in}(y)|$ we need appropriate estimate on a length 
of a characteristic connecting $y$ with the boundary and here the assumptions \eqref{flat} 
and (A1) come into play. The required result is  
\begin{lem} \label{L6}
Let $b(t)=(t,\uxd(-t))$ and let $p(t)$ satisfy 
\begin{equation} \label{eq:p}
\dot p(t)=(1+U^1(p(t)),U^2(p(t))), \quad p(t_0)=x_2^0.
\end{equation}
where $x_2^0$ is a given, small positive number, $t_0$ is sufficiently small that 
$$
p(s)>b(s) \quad {\rm for} \quad s \leq t_0. 
$$
and $U$ satisfies the assumptions of Proposition \ref{lem_trans}. 
Let $t_*$ denote the first coordinate of the first intersection of $p(t)$ and $b(t)$ (see Fig. 5). 
Then  
\begin{equation*}
t_*-t_0 \leq C (x_2^0-b(t_0))^\delta,
\end{equation*} 
where $\delta$ is from \eqref{flat_inv} and $C=C(U,\Omega)$.
\end{lem}
\begin{rem} We will apply the above lemma with $t=-x_1$ in some neighbourhood of a singularity point 
$(0,0)$ in the proof of Lemma \ref{L7} which will follow, therefore we consider a 'backward' characteristic
of the equation \eqref{trans}. 
\end{rem}
\emph{Proof.}
Assume first that $t_0=0$. 
Let us denote 
$$
A(t)=\sqrt{1+|\uxd'(t)|^2}. %\quad A_2(t,z)=\sqrt{(1+U_1(t,z))^2+(U_2(t,z))^2}.
$$
Then the unit normal to $b(t)$ is 
$n_b(t)={A(t)}^{-1}(-\uxd'(t),1).$
Let $z(t)$ denote the distance between $b(t)$ and $p(t)$ measured along $n_b(t)$ (see Fig. 5), that is 
\begin{equation*}
p(t)=b(t)+n_b(t)z(t).
\end{equation*}
Differentiating this identity in $t$ and multiplying the resulting equation by $n_b(t)$ we get
\begin{equation} \label{zdot:1}
\dot{p}(t) \cdot n_b(t)=\dot{z}(t),
\end{equation} 
%On the other hand, 
%$$
%\dot{p}(t)=\frac{1}{A_2}(1+U_1(t,z),U_2(t,z))=\frac{1}{A_2}(1,0)+\frac{1}{A_2}(U_1(t,z),U_2(t,z))
%$$
Which by \eqref{eq:p} rewrites
\begin{equation} \label{zdot:2}
\dot{z}(t)=U \cdot n_b(t,z)-\frac{\underline{x_2}'(t)}{A(t)}.
\end{equation}
We have 
$$
|U\cdot n_b(t,z)| \leq |U\cdot n_b(t,0)|+\|\nabla U\|_{L_\infty}z, 
$$
therefore by \eqref{zdot:2}
\begin{equation} \label{zdot:3}
\dot{z}(t) \leq |U\cdot n_b(t,0)|-\frac{\uxd'(t)}{A(t)}+\|\nabla U\|_{L_\infty}z(t).
\end{equation}
Now, due to \eqref{infl:small} and \eqref{def:lambda} we have 
\begin{equation} \label{zdot:4}
|U \cdot n_b(t,0)| \leq \frac{\tilde \kappa |\uxd'(t)|}{\sqrt{1+\lambda^2}} < \tilde \kappa\frac{|\uxd'(t)|}{A(t)}
\end{equation}
therefore  
\begin{equation} \label{zdot:4}
\dot{z}(t) \leq - (1- \tilde \kappa) \frac{\uxd'(t)}{A(t)}+E_U z(t) \leq -\theta +E_U z(t)  
\end{equation}
with $\theta=\frac{1-\tilde \kappa}{\sqrt{1+\lambda^2}}$ %, where $\lambda$ is from \eqref{def:lambda}, 
and $E_U=\|\nabla U\|_{L_\infty}$.
Now $t_*$ must satisfy 
\begin{equation} \label{zdot:5}
\int_0^{t^*}\dot z(t)dt = -x_2^0,
% \end{equation}
\mbox{ \ therefore \ } 
% \begin{equation} \label{zdot:5}
x_2^0=\theta \uxd(t_*)-E_U\underbrace{\int_0^{t_*}z(t)dt}_{Z(t_*)}.  
\end{equation}
In order to simplify the above condition observe that 
$$z(t)\leq x_2^0 e^{E_U t} \leq C x_2^0, \quad  p(t)-b(t) \leq C z(t)$$,
and therefore   
\begin{equation} \label{eZ}
\displaystyle Z(t_*)\leq t \,\sup_{s \leq t}[p(t)-b(t)] \leq C_2 t_* x_2^0.
\end{equation}
Moreover, by \eqref{flat} we have 
\begin{equation} \label{flat:3}
\uxd(t_*)\geq C_1 t_*^N.
\end{equation}
In view of \eqref{eZ} and \eqref{flat:3}, if we define $\tilde t_*$ by 
\begin{equation*}
P_{x_2^0}(\tilde t_*):=C_1\theta \tilde t_*^N - C_2 E_U x_2^0 \tilde t_* - x_2^0 = 0, 
\end{equation*}
then, because of \eqref{zdot:5},\eqref{eZ} and \eqref{flat:3},
\begin{equation*} \label{eT}
t_* \leq \tilde t_*
\end{equation*} 
We claim that there exists $1 \leq C_3=C_3(C_1,C_2,\Omega)$ s.t. if
$t_1=C_3 (x_2^0)^{1/N}$ 
then 
% \begin{equation}
$P(t_1) \geq 0.$
% \end{equation}
Indeed, we have
$$ 
P(t_1)= x_2^0 [ C_1\theta C_3-1-C_2 E_U (C_3 x_2^0)^{1/N}]\geq 
C_3 (C_1-C_2 E_U (x_2^0)^{1/N})-1 >0 
$$
for sufficiently large $C_3=C_3(C_1,C_2,\Omega)$. 
Therefore, 
\begin{equation*}
\tilde t_* \leq t_1 \leq C_3 (x_2^0)^{1/N}, 
\end{equation*}
which together with \eqref{eT} completes the proof for $t_0=0$ as $\delta=\frac{1}{N}$. 
For $t_0>0$ the reasoning is the same modulo minor adjustments,
replacing $x_2^0$ with $-z(t_0)$ in \eqref{zdot:5}.  
\begin{figure} 
\begin{center}
\includegraphics[width=0.8\textwidth]{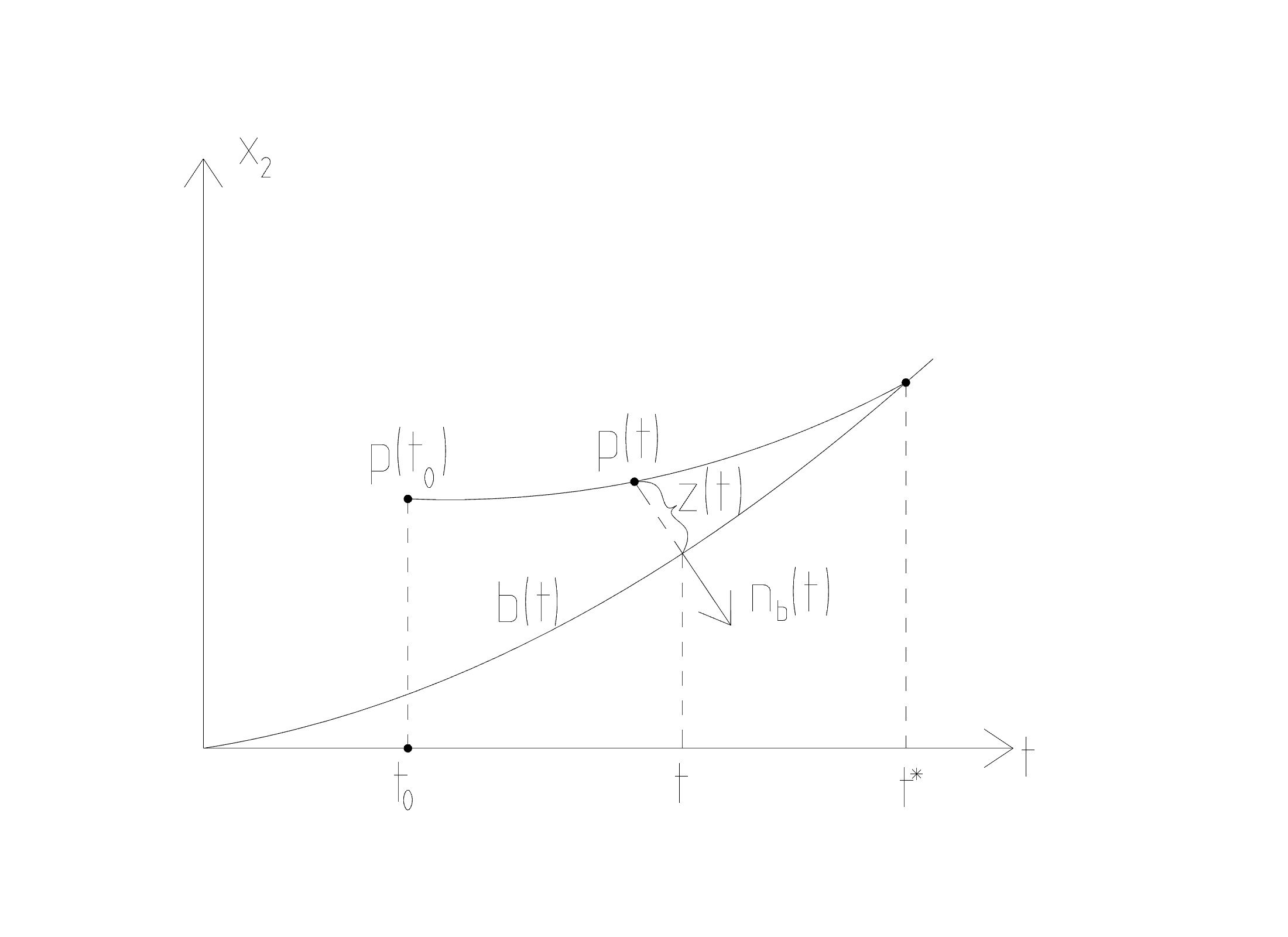}
\caption{Explanation of notation used in Lemma 6}
\end{center}
\end{figure}
\noindent
We are now in a position to prove the estimate on the boundary term on the lhs of \eqref{est_trans_2}.  
\begin{lem} \label{L7}
Let the assuptions of Proposition \ref{lem_trans} hold. Then
\begin{align} \label{est_bd1}
&\int_{\Omega} dy \int_{\Gamma_{in}} \frac{|w(x)-w(y)|^p (n_1+U \cdot n)(x)}{\phi_{\epsilon}(x,y)} dS(x) \leq\nonumber\\
&\leq C [\|H\|_{L_\infty(\Omega)}+\|w\|_{L_p(\Omega)}+\|w\|_{L_\infty(\Omega)}+\|w_{in}\|_{W^s_p(\Gamma_{in})}
+\|w_{in}\|_{W^r_p(\Gamma^s_{in})})]^p.
\end{align}
\end{lem}
\emph{Proof.} 
%Recalling the definition of sets $\Gamma_{in}^r$ and $\Gamma_{in}^s$ introduced before Lemma \ref{L5}, 
%we have 
%$$
%\left|\frac{U \cdot n}{n_1}\right|\leq \|U\|_{\infty}\sqrt{1+[\ux'(\eta)]^2} \quad {\rm on} \quad \Gamma_{in}^r
%$$
%and 
%$$
%\left|\frac{U \cdot n}{n_1}\right|\leq \kappa \frac{\sqrt{1+\ux'(x_2)}}{|\ux'(x_2)|} 
%\leq \kappa \frac{\sqrt{1+\ux'(\eta)}}{|\ux'(\eta)|},
%$$
First of all, due to \eqref{normal} and (A1) we have
$$
|\frac{U\cdot n}{n_1}| \leq \kappa,
$$ 
which gives
\begin{equation} \label{infl}
|n_1+U \cdot n| \leq (1+\kappa)n_1.
\end{equation} 
Therefore it is enough to show \eqref{est_bd1} for $U \cdot n=0$.
As 
\begin{equation} \label{normal_0}
dS(x_2)|_{\Gamma_{out}}=\sqrt{1+[\overline{x_1}'(x_2)]^2}dx_2 \quad \textrm{and} \quad
dS(x_2)|_{\Gamma_{in}}=\sqrt{1+[\underline{x_1}'(x_2)]^2}dx_2,
\end{equation}
\eqref{normal} yields
\begin{equation} \label{normal_1}
dS(x_2)|_{\Gamma_{in}} = -\frac{dx_2}{n_1}, \qquad n_2 dS(x_2)|_{\Gamma_{in}}=\underline{x_1}'(x_2)dx_2.
\end{equation}
By \eqref{normal_1}, assuming $U \cdot n=0$ the LHS of \eqref{est_bd1} can be rewritten as
\begin{equation} \label{w1}
\int_{\Omega} \int_{\Gamma_{in}} \frac{|w(x)-w(y)|^p n_1(x)}{\phi_{\epsilon}(x,y)} dS(x)dy=\int_{\Omega} dy \underbrace{\int_0^b 
- \frac{|w(\underline{x_1}(x_2),x_2)-w(y)|^p}{\phi_{\epsilon}((\underline{x_1}(x_2),x_2),y)} \,dx_2}_{I_\epsilon(y)}, 
\end{equation}
and it remains to show the estimate \eqref{est_bd1} for \eqref{w1}. 
When we pass with $\epsilon \to 0$, we can expect some problems only when $y$ is close to $\Gamma_{in}$. 
Therefore we write \eqref{w1} as 
$$
\int_{\Omega} I_\ep(y) dy  = 
\int_{\Omega^r_{in}} I_\ep(y) dy + \int_{\Omega^s_{in}} I_\ep(y) dy + \int_{\Omega \setminus \Omega_{in}} I_\ep(y) dy  =:
I_1^r + I_1^s + I_2,
$$
where the sets $\Omega^r_{in}$, $\Omega^s_{in}$ and $\Omega_{in}$ are defined in \eqref{def:omega_in}.
The set $\Omega\setminus \Omega_{in}$ consist of points at the distance from $\Gamma_{in}$ bounded 
from below by $\eta$. Therefore the $I_2$ integral is straightforward as we have $|(\underline{x_1}(x_2),x_2)-y|>\eta$ for $y \in \Omega \setminus \Omega_{in}$ and so
\begin{equation} \label{i2}
I_2 \leq C(\eta) [ \|w_{in}\|_{L_p(\Gamma_{in})} + \|w\|_{L_p} ]^p. 
\end{equation}
The estimates for $I^r_2$ and $I^s_1$ are more involved, we show them in details. \\
\noindent
{\bf Estimate of $I^r_1$}. 
Let $\underline{x}(x_2)=(\underline{x_1}(x_2),x_2) \in \Gamma_{in}$.
We have
\begin{equation} \label{br}
\frac{|w(\underline{x}(x_2))-w(y)|^p}{|\underline{x}(x_2)-y|^{2+sp}} \leq C(p) \left[
\frac{|w(\underline{x}(x_2))-w(\uyy(y))|^p}{|\underline{x}(x_2)-y|^{2+sp}}
+ \frac{|w(\uyy(y))-w(y)|^p}{|\underline{x}(x_2)-y|^{2+sp}} \right]. 
\end{equation}
%\frac{|w(\underline{x})-w(\underline{y})|^p}{|\underline{x}-\underline{y}|^{2+sp}}
%+ \frac{|w(\underline{y})-w(y)|^p}{|\underline{x}-y|^{2+sp}}, 
%where in the second inequality we used either \eqref{br1} or \eqref{br2}.
For the first term on the rhs of \eqref{br} we have
\begin{align} \label{sing:3} 
&\int_{\Omega^r_{in}}dy \int_0^b \frac{|w(\underline{x}(x_2))-w(\uyy(y))|^p}{|\underline{x}(x_2)-y|^{2+sp}} \,dx_2 = \nonumber\\
&=\int_{2\eta}^{b-2\eta} dx_2 \left[ \int_{V_\eta(x_2)} \frac{|w(\underline{x}(x_2))-w(\uyy(y))|^p}{|\underline{x}(x_2)-y|^{2+sp}}\,dy  
+ \int_{\Omega^r_{in} \setminus V_\eta(x_2)} \frac{|w(\underline{x}(x_2))-w(\uyy(y))|^p}{|\underline{x}(x_2)-y|^{2+sp}}\,dy \right],
\end{align}
where
\begin{equation} \label{def:v_eta}
V_\eta(x_2)=\{ y \in \Omega^r_{in}: \; x_2-\eta < y_2 < x_2+\eta \},
\end{equation} 
which implies that in the  second integral $|\uxx(x_2)-y|$ is bounded from below by $\eta$, therefore
%Moreover, since we are away from the 
%singularity points, we have 
%Therefore we obtain
\begin{align*}
&\int_0^b dx_2\int_{\Omega^r_{in} \setminus V_\eta(x_2)} \frac{|w(\underline{x}(x_2))-w(\uyy(y))|^p}{|\underline{x}(x_2)-y|^{2+sp}}\,dy \nonumber\\ 
&\leq C \left[\int_0^b |w(\uxx(x_2))|^p dx_2 + \int_0^b |w(\uyy(y))|^p dy_2\right] \leq C \|w_{in}\|_{L_p(\Gamma_{in})},  
\end{align*}
where in the last passage we have used an obvious inequality 
\begin{equation} \label{sing:5}  
dy_2 \leq dS(y_2).
\end{equation} 
In order to treat the first integral in \eqref{sing:3} we apply \eqref{sing:1}-\eqref{sing:1b}.
Due to \eqref{sing:1b} we can assume w.l.o.g. for this estimate that $\gamma_y$ is a straight line $y_2 \equiv {\rm const}$,
therefore $\uyy(y)=\underline{y}(y)$. Then   
by the definitions of $V_\eta(x_2)$ and $\Omega^r_{in}$ (\eqref{def:v_eta} and \eqref{def:omega_in}, respectively)   
\begin{align*}    
&\int_{V_\eta(x_2)} \frac{|w(\underline{x}(x_2))-w(\uyy(y))|^p}{|\underline{x}(x_2)-y|^{2+sp}}\,dy=
\int_{x_2-\eta}^{x_2+\eta}\int_{\uy(y_2)}^{\uy(y_2)+\eta}\frac{|w(\uxx(x_2))-w(\uyy(y_2))|^p}{|\uxx(x_2)-y|^{2+sp}}dy_1dy_2 \nonumber\\
&\leq \int_{x_2-\eta}^{x_2+\eta}|w(\uxx(x_2))-w(\uyy(y_2))|^pdy_2\int_{\uy(y_2)}^{\uy(y_2)+\eta}(y_1-\uy(y)+|\uxx(x_2)-\uyy(y_2)|)^{-2-sp}dy_1\nonumber\\ 
&\leq \int_{x_2-\eta}^{x_2+\eta}\frac{|w(\uxx(x_2))-w(\uyy(y_2))|^p}{|\uxx(x_2)-\uyy(y_2)|^{1+sp}}dy_2,
\end{align*}
where in the first inequality we have used \eqref{sing:1}.
Therefore, as we have additionally \eqref{sing:5},
\begin{equation} \label{i1r_1}
\int_{2\eta}^{b-2\eta} dx_2\int_{V_\eta(x_2)} \frac{|w(\underline{x}(x_2))-w(\uyy(y))|^p}{|\underline{x}(x_2)-y|^{2+sp}}\,dy\leq 
C\|w\|_{W^s_p(\Gamma_{in})}.
\end{equation}
Notice that in this part of the estimate the norm of boundary data appears naturally.

In the second term on the r.h.s. of \eqref{br} we express $w(y)-w(\uyy(y))$ by an integral along $\gamma_y$ 
using the equation \eqref{trans}.
Applying Jensen inequality we get
\begin{align} \label{i1r_2}
&\int_{\Omega^r_{in}} dy \int_0^b \frac{|w(\uyy(y))-w(y)|^p}{|\underline{x}(x_2)-y|^{2+sp}} \,dx_2   \leq \int_{\Omega_{in}^r} dy \int_0^b \frac{ |\int_{\gamma_y} (H-w) dl_{\gamma_y}|^p }{ |\underline{x}(x_2)-y|^{2+sp} }\,dx_2 \leq\nonumber
\\[3pt]
&\leq(\|H\|_{L_\infty}+\|w\|_{L_\infty})^p \int_{\Omega_{in}^r} dy \int_0^b \frac{|\gamma_y|^p}{|\underline{x}(x_2)-y|^{2+sp}} \,dx_2 \leq \nonumber
\\[3pt]
&\leq C(\|H\|_{L_\infty}+\|w\|_{L_\infty})^p \int_{\Omega_{in}^r} dy \int_0^b \frac{|\uyy(y) - y|^p}{|\underline{x}(x_2)-y|^{2+sp}} \,dx_2 \leq \nonumber
\\[3pt]
&\leq C(\|H\|_{L_\infty}+\|w\|_{L_\infty})^p \int_0^b dx_2 \int_{\Omega_{in}^r} \frac{|\underline{x}(x_2) - y|^p}{|\underline{x}(x_2)-y|^{2+sp}} \,dy, 
\end{align}
where we have used \eqref{sing:1} and \eqref{sing:1b}.
The last integral is finite for $s<1$, we remember $\Omega^r_{in}$ is two dimensional.
Combining \eqref{br}, \eqref{i1r_1} and \eqref{i1r_2} we get
\begin{equation} \label{i1r}
I_1^r \leq C (\|H\|_{L_\infty(\Omega)}+\|w\|_{L_\infty(\Omega)} + \|w_{in}\|_{W^s_p(\Gamma_{in})})^p.
\end{equation}

{\bf Estimate of $I^s_1$}.
Around the singularity points we have to proceed more carefully. 
If we consider again separately the cases \eqref{cases},
then in the case {\bf (a)} \eqref{sing:1} still holds with $C=2$
so we can repeat directly the previous approach and  
it remains to consider the cases {\bf (b)} and {\bf (c)}. 
The key difference is that, as we already explained, 
$|\ux'(x_2)|$ is unbounded.
\begin{figure} 
\begin{center}
\includegraphics[width=0.8\textwidth]{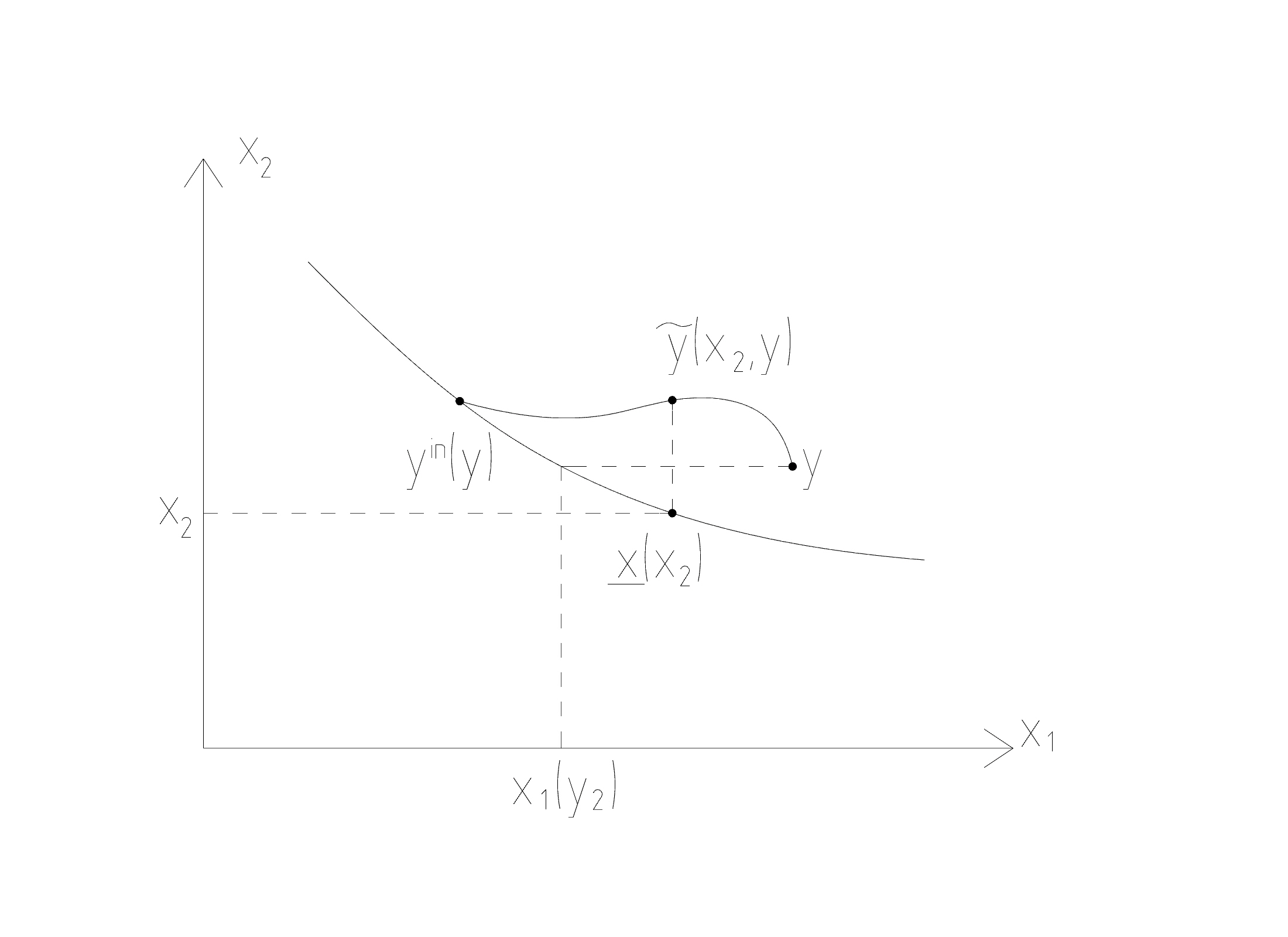}
\caption{Illustration of notation. Case {\bf (b)}: intersection point $\tilde y$ well defined.}
\end{center}
\end{figure}
\begin{figure} 
\begin{center}
\includegraphics[width=0.8\textwidth]{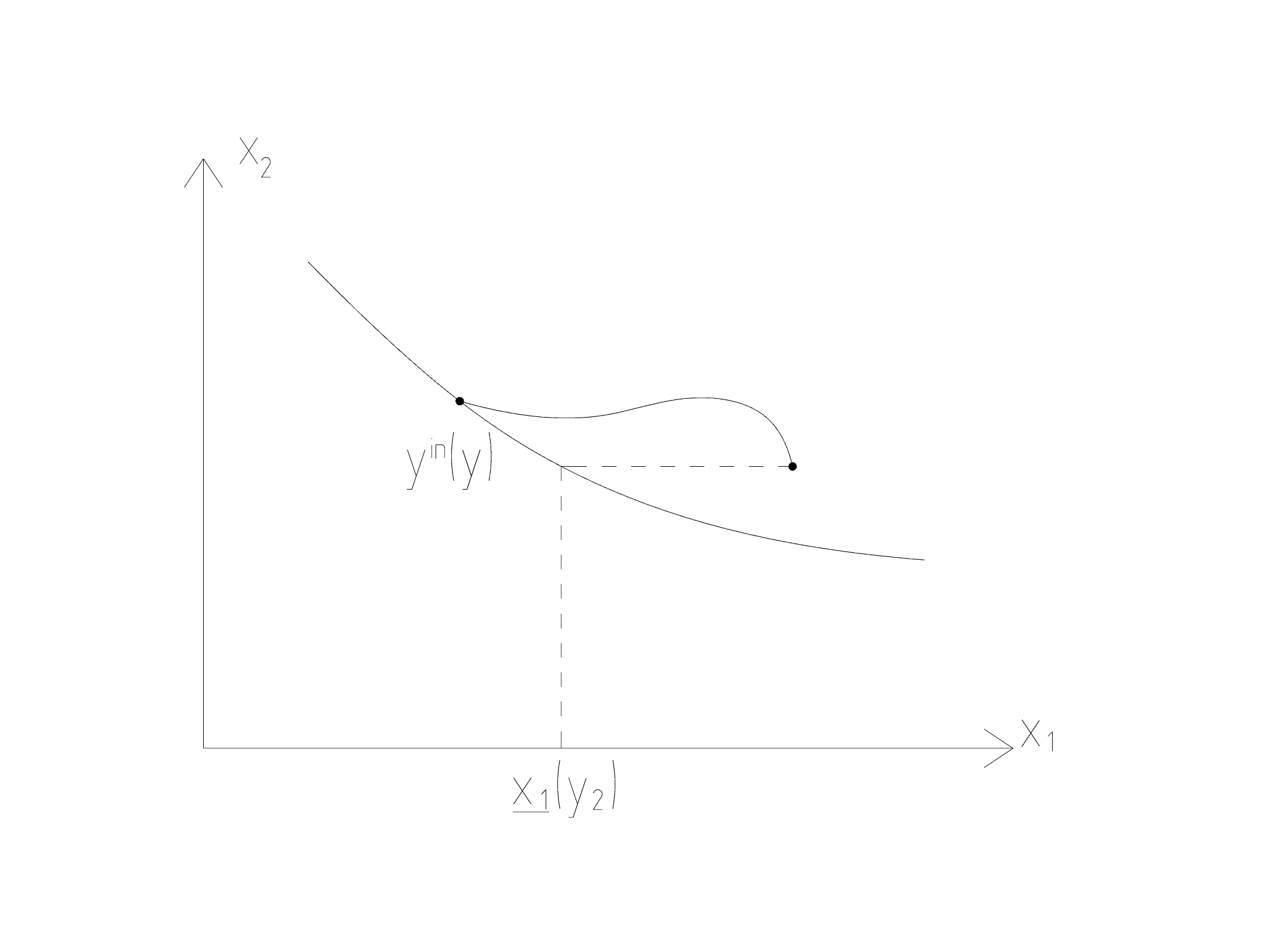}
\caption{Case {\bf (a)}: no intersection point $\tilde y$.}
\end{center}
\end{figure}
Therefore we cannot repeat directly \eqref{i1r_1} and \eqref{i1r_2}. 
We have to control $|\underline{x}-y^{in}(y)|$ and here we use the assumption \eqref{flat_inv}.
Recall the notions of $\gamma_y$ and $\tilde y$ introduced in \eqref{def:gamma} and after \eqref{cases},
respectively, see Fig. 3.  
Let us denote by $\gamma_{\tilde y}$ the part of the $\gamma_y$ connecting $\Gamma_{in}$ with $\tilde y(x_2,y)$. 
We modify \eqref{br} adding the point $\tilde y(x_2,y)$:
\begin{align} \label{bs}
\frac{|w(\underline{x}(x_2))-w(y)|^p}{|\underline{x}(x_2)-y|^{2+sp}} 
\leq & C(p) \left[  \frac{|w(\underline{x}(x_2))-w(\uyy(y))|^p}{|\underline{x}(x_2)-y|^{2+sp}}
+ \frac{|w(\uyy(y))-w(\tilde y(x_2,y))|^p}{|\underline{x}(x_2)-y|^{2+sp}} \right.  \nonumber\\
& \left. + \frac{|w(\tilde y(x_2,y))-w(y)|^p}{|\underline{x}(x_2)-y|^{2+sp}} \right]. 
\end{align}
%$$
%\frac{|w(\underline{x})-w(\tilde y)|^p}{|\underline{x}-\tilde y|^{2+sp}}
%+ \frac{|w(\tilde y)-w(y)|^p}{|\tilde y-y|^{2+sp}}.
%$$
We start with the 2nd and 3rd term on the rhs of \eqref{bs}.
Similarly as before, in both terms we replace the value of $w$ with an integral along $\gamma_y$.
The last term is analogous to $I_1^r$ and we get 
\begin{equation} \label{is1}
\int_{\Omega^s_{in}}\,dy\int_0^b \frac{|w(\tilde y(x_2,y))-w(y)|^p}{|\underline{x}(x_2)-y|^{2+sp}} \,dx_2 \leq 
C ( \|H\|_{L_{\infty}(\Omega)}+\|w\|_{L_{\infty}(\Omega)} )^p \quad \textrm{for} \quad s<1 .
\end{equation}
In the second term we have 
\begin{align*}
|w(\tilde y(x_2,y)) - w(\uyy(y))| &= |\int_{\gamma_{\tilde y}} (H-w) dl_{\gamma_{\tilde y}}|\\
&\leq (\|H\|_{L_{\infty}}+\|w\|_{L_{\infty}})\int_{y^{in}_1(y)}^{\ux(x_2)}\sqrt{[1+U_1(\gamma_y(s))]^2+(U_2(\gamma_y(s)))^2}{\rm ds}\\
&\leq C(\|H\|_{L_{\infty}}+\|w\|_{L_{\infty}}) |y^{in}_1(y)-\ux(x_2)|.
\end{align*}
Now we  estimate the latter term. This is the most subtle part of the whole 
proof which has been carried out in Lemma \ref{L6}. 
Applying this result with $t_0=-\ux(x_2)$ we have 
\begin{equation} \label{est:char}
|y^{in}_1(y)-\ux(x_2)| \leq C |x_2 - \tilde y_2(x_2,y)|^\delta,
\end{equation}
where $\delta$ is from \eqref{flat_inv}.
Next,
$$
|x_2-\tilde y_2(x_2,y)|\leq |x_2-y_2|+\|U_2\|_{L_\infty(\Omega)}|x_1-y_1|\leq C |\uxx-y|,
$$
therefore
\begin{align*}
&\int_{\Omega^s_{in}} dy \int_0^b \frac{|w(\uyy(y))-w(\tilde y(x_2,y))|^p}{|\underline{x}(x_2)-y|^{2+sp}} \leq
\int_{\Omega^s_{in}} dy \int_0^b \frac{|\int_{\gamma_{\tilde y}}(H-w)dl_{\gamma_{\tilde y}}|^p}{|\underline{x}(x_2)-y|^{2+sp}}\,dx_2 \leq 
\\[3pt]
&\leq C (\|H\|_{L_{\infty}}+\|w\|_{L_{\infty}})^p \int_{\Omega^s_{in}} dy \int_0^b \frac{|x_2-\tilde y_2(x_2,y)|^{\delta p}}{|\underline{x}(x_2)-y|^{2+sp}}\,dx_2 \leq
\\[3pt]
& \leq C (\|H\|_{L_{\infty}}+\|w\|_{L_{\infty}})^p \int_{\Omega^s_{in}} dy \int_0^b |\uxx(x_2)-y|^{\delta p-2-sp} \,dx_2,
\end{align*}
where in the second passage we used \eqref{est:char}.
The last integral is finite for $s<\delta$ and we conclude
\begin{equation} \label{is2}
\int_{\Omega^s_{in}}\,dy\int_0^b \frac{|w(\uyy(y))-w(\tilde y(x_2,y))|^p}{|\underline{x}(x_2)-y|^{2+sp}} \,dx_2 \leq 
C ( \|H\|_{L_{\infty}}+\|w\|_{L_{\infty}} )^p \quad \textrm{for} \quad s<\delta .
\end{equation}
It remains to estimate the first term on the rhs of \eqref{bs}.
For this purpose we  assume $\underline{x} \neq \uyy(y)$ since otherwise 
this term vanishes. 
%
%Moreover, around the singularity points we have
%\begin{equation} \label{sing:18}
%\frac{dx_2}{dS} \leq C \frac{|x_2-y_2|}{|\uxx(x_2)-\underline{y}(y)|}.
%\end{equation}
%Moreover, be Lemma \ref{L6} we have \\
%
%\smallskip
%\smallskip
%$$
%|\uxx(x_2)-\underline{y}(y)| \leq C |\uxx(x_2)-\underline{y}(y)|,
%$$
%therefore from \eqref{sing:18}
%\begin{equation} \label{sing:18}
%\frac{dx_2}{dS} \leq C \frac{|x_2-y_2|}{|\uxx(x_2)-\uyy(y)|}\leq C \frac{|\uxx(x_2)-y|}{|\uxx(x_2)-\uyy(y)|}
%\end{equation}
%
%\begin{equation} \label{ds_sing}
%\frac{dx_2}{dS} \leq C \frac{|x_2-y_2|}{|\uxx(x_2)-\uyy(y)|} \leq C \frac{|\uxx(x_2)-y|}{|\uxx(x_2)-\uyy(y)|}.  
%\end{equation}
%
By Fubini's theorem we have $dy=dS(y_2)dy_1$, therefore we can write
\begin{align} \label{sing:19}
&\int_{\Omega^s_{in}} dy \int_0^b \frac{|w(\underline{x}(x_2))-w(\uyy(y))|^p}{|\underline{x}(x_2)-y|^{2+sp}}dx_2\nonumber\\
&\leq C \int_{\Gamma^s_{in}} dS(y_2) \int_0^b dx_2 \int_{\underline{y_1}(y_2)}^{\overline{y_1}(y_2)}
\frac{|w(\underline{x}(x_2))-w(\underline{y}(y_2)))|^p}{[|\ux(x_2)-y_1|+|x_2-y_2|]^{2+sp}}dy_1\nonumber\\
&\leq C \int_{\Gamma^s_{in}} dS(y_2) \int_0^b dx_2 
\frac{|w(\underline{x}(x_2))-w(\underline{y}(y_2))|^p}{|x_2-y_2|^{1+sp}}.
\end{align}
Introducing $h=x_2-y_2$ and using the imbedding \eqref{imbed:2} we get 
\begin{align} \label{sing:20} 
&\int_{\Gamma^s_{in}} dS(y_2) \int_0^b dx_2 \frac{|w(\underline{x}(x_2))-w(\underline{y}(y_2))|^p}{|x_2-y_2|^{1+sp}}   \nonumber\\
&\leq C \int_0^{h_0} \frac{dh}{|h|^{1+sp}} \int_{\Gamma^{s,h}_{in}|w(\uxx(y_2+h))-w(\underline{y}(y_2))|^p dS(y_2)}\\
&\leq C \|w\|_{B^r_{p,\infty}(\Gamma^s_{in})}\int_0^{h_0} h^{-1-sp} {\rm sup}_{y_2 \in (0,\eta)}|\uxx(y_2+h)-\underline{y}(y_2)|^{rp}dh \nonumber\\
&\leq C \|w\|_{W^r_p(\Gamma^s_{in})}\int_0^{h_0} h^{-1-sp} {\rm sup}_{y_2 \in (0,\eta)}|\uxx(y_2+h)-\underline{y}(y_2)|^{rp}dh ,\nonumber
\end{align} 
where 
$$
h_0={\rm sup}\{ x_2: \; (\ux(x_2),x_2) \in \Gamma^s_{in}\} \quad {\rm and} \quad 
{\Gamma^{s,h}_{in}}=\{ \uxx(x_2): \; \uxx(x_2+h)\in \Gamma^s_{in}  \}. 
$$
In the second passage in \eqref{sing:20} we have used \eqref{norm:brp} and in the last one \eqref{imbed:2}.
Now we apply Lemma \ref{L6} with 
$$
t_0=-x_1, \quad p(t_0)=(\uyy_1(y_2),y_2+h), \quad p(t)=\gamma_{p(t_0)}(-t),
$$ 
where $\gamma_y(\cdot)$ is a characteristic passing through $y$ defined in \eqref{def:gamma}.
Then Lemma \ref{L6} implies 
\begin{equation} \label{sing:21}
|\uxx(y_2+h)-\underline{y}(y_2)| \leq C|\uxx(y_2+h)-p(t_0)| \leq C|h|^\delta,
\end{equation}
which together with \eqref{sing:20} gives 
\begin{equation} \label{is3}
\int_{\Gamma^s_{in}} dS(y_2) \int_0^b dx_2 \frac{|w(\underline{x}(x_2))-w(\underline{y}(y_2))|^p}{|x_2-y_2|^{1+sp}}
\leq C \|w\|_{W^r_p(\Gamma^s_{in})} \int_0^{h_0} h^{\delta rp-1-sp}.
\end{equation}
The last integral is finite provided $$\delta r>s$$,
which implies the second relation in \eqref{c:sp}. 
%
%In the above series of inequalities first we used twice \eqref{ds_sing} and then \eqref{flat_inv_1}.
%As obviously $2+\frac{sp}{\delta}>1+sp$, we need higher regularity of $w_{in}$ on $\Gamma^s_{in}$. 
%Assuming $W \in W^r_p(\Gamma^s_{in})$ we need
%$$
%1+\frac{sp}{\delta}<rp, 
%$$
%In particular, for $r=1$ we get \eqref{c:rhoin2}.
%Under these assumptions we get from \eqref{sing:19} and \eqref{sing:20}
%\begin{equation} \label{is3}
%\int_{\Omega^s_{in}} dy \int_0^b \frac{|w(\underline{x}(x_2))-w(\uyy(y))|^p}{|\underline{x}(x_2)-y|^{2+sp}}dx_2 \leq
%\|w_{in}\|_{W^r_p(\Gamma_s)}^p. 
%\end{equation}
Combining \eqref{is1}, \eqref{is2} and \eqref{is3} we conclude
\begin{equation} \label{i1s}
I^s_1 \leq C (\|H\|_{L_{\infty}}+\|w\|_{L_{\infty}}+\|w_{in}\|_{W^r_p(\Gamma_s)})^p.  
\end{equation}
Putting together \eqref{i2},\eqref{i1r} and \eqref{i1s} we get \eqref{est_bd1} with $U \cdot n=0$,
which completes the proof due to \eqref{infl}. 

\qed

\emph{Proof of Proposition \ref{lem_trans}}. Now we are ready to close the estimate \eqref{est_trans}. 
Combining \eqref{est_trans_2} with \eqref{est_bd1} we obtain 
\begin{align} \label{est_trans_0}
\|w\|_{W^s_p}^p \leq & C(s,p,\Omega) [ \|H\|_{W^s_p} \|w\|_{W^s_p}^{p-1} + \|w\|_{L_\infty}^p + \|w\|_{L_p}^p \nonumber\\ 
& + \|w_{in}\|_{W^s_p(\Gamma_{in})} + \|w_{in}\|_{W^r_p(\Gamma_{s})} + \|\nabla U\|_{L_{\infty}} \|w\|_{W^s_p}^p]
\end{align}
and applying the interpolation inequality \eqref{int} we conclude \eqref{est_trans}.
\qed
\begin{rem} \label{rem_indep} 
Notice that from \eqref{est_trans_0} it follows that the constant in \eqref{est_trans} 
is of a form $C=\frac{C(s,p,\Omega)}{1-\|\nabla U\|_{L_\infty}}$, therefore under the assumption of smallness of $\|\nabla U\|_{L_\infty}$ we can assume that 
it does not depend on $\|\nabla U\|_{L_\infty}$.
\end{rem}
\begin{rem}
The proof of Proposition \ref{lem_trans} is quite technical, therefore it may be unclear for the reader how 
it can be extended to more general domain. To clarify it, we shall emphasize that the most subtle 
part is the estimate 
of $\frac{|w(\uxx)-w(y)|^p}{|\uxx-y|^{2+sp}}$ where $\uxx$ is on the boundary and $y$ is 
inside $\Omega$, both close to singularity. 
The most delicate point is to estimate the distance between the point where a characteristic crossing $y$ 
intersects the boundary and $\uxx$. This was carried out in Lemma \ref{L6}. 
On the other hand, larger distance between $\uxx$ and $y$ does not 
entail additional problems. Therefore the proof presented here can be extended to the general 
case of multiple singularity points. We have to consider again separately the 'regular'
part of the boundary where $\uxd'(x_1)$ is bounded by a given constant and therefore 
we get the estimates \eqref{i2}, \eqref{i1r}, and neighbourhoods of 
each singularity point where we repeat the estimate \eqref{i1s}.   
\end{rem}

\qed

\subsection{Linear estimate in $W^{1+s}_p \times W^s_p$ and solution of the linear system}
In this section we solve linear problem \eqref{system_lin}
with rhs of regularity determined in \eqref{reg_rhs}. First we show the \emph{a priori} estimate.
\begin{prop} \label{prop_apriori}
Let $(u,w)$ solve the system \eqref{system_lin} with the rhs of regularity determined in \eqref{reg_rhs} 
with $s,p$ and the boundary data satisfying the assumptions of Theorem 1. Then we have
\begin{equation} \label{est_lin}
 \|u\|_{W^{1+s}_p(\Omega)} + \|w\|_{W^s_p(\Omega)} \leq \\
C \, [\|F\|_{L_p(\Omega)} + \|G\|_{W^s_p(\Omega)}+\|w_{in}\|_{W^s_p(\Gamma_{in}) \cap W^r_p(\Gamma_s)}
+ \|B\|_{W^{1-1/p}_p(\Gamma)} ]. 
\end{equation}
\end{prop}
%\begin{rem}
%Notice the for the estimate we do not need the assumption (A3).
%\end{rem}
\begin{rem}
The main difficulty in the proof of \eqref{est_lin} lies in the steady transport equation which 
has been dealt with in Proposition \ref{lem_trans}. The rest of the proof uses classical results from the 
theory of elliptic problems and can be divided in several steps.
\end{rem}
\emph{Proof of Proposition \ref{prop_apriori}}. 
{\bf Step I.} 
Let us take the rotation of $(\ref{system_lin})_1$. Then we get the following system
\begin{equation} \label{AP2}
 \begin{array}{lcr}
  -\mu \Delta \rot u = \rot (F-\d_{x_1}u) & \mbox{in} & \Omega,\\[3pt]
\rot u = (2\chi - f/\nu ) u\cdot \tau + \frac{B}{\mu} & \mbox{ at } \Gamma.
 \end{array}
\end{equation}
In order to show the boundary relation for the vorticity we differentiate 
\eqref{system_lin}$_3$ in the tangential direction and apply \eqref{system_lin}$_4$.
The details are shown in \cite{MR} in the case 
$u \cdot n = d, \; B=0$. A minor modification for our case $u \cdot n = 0, B \neq 0$
yields \eqref{AP2}.
%Here we use the form of $F$ specified in \eqref{reg_rhs}. 
%Since $F=\nabla F_1+f_2$, 
We find the following estimate for the solution to (\ref{AP2}):
\begin{equation} \label{AP3}
 \|\rot u\|_{W^1_p(\Omega)} \leq C( \|F\|_{L_p(\Omega)} + \|u,B\|_{W^{1-1/p}_p(\Gamma)} + \|\d_{x_1}u\|_{L_p(\Omega)} ).
\end{equation}
The solvability of (\ref{AP2}) belongs to the classical theory of elliptic linear problems, hence we omit the details of this construction.

\bigskip

{\bf Step II.} Consider the Helmholtz decomposition of $u$:
\begin{equation} \label{Helm}
u = \nabla \phi + \nabla^\bot A.
\end{equation} 
We have $0=n \cdot \nabla^\bot A=\tau \cdot \nabla A$. Therefore $A$ satifies 
\begin{equation} \label{eq:A}
\Delta A = \rot u, \quad A|_{\Gamma}={\rm const}.
\end{equation}
As $\partial \Omega \in C^3$, 
from \eqref{eq:A} we obtain 
\begin{equation} \label{AP4}
\|A\|_{W^3_p(\Omega)} \leq C \|\rot u\|_{W^1_p(\Omega)}.
\end{equation}
Substituting \eqref{Helm} to $(\ref{system_lin})_1$ we get 
\begin{equation*}
\nabla(-(2\mu+\nu)\div u + Kw)=F-\de_{x_1}u+\mu \Delta \nabla^\bot A + (\mu+\nu)\nabla \div \nabla^\bot A=:\bar F.
\end{equation*}
By \eqref{AP3} and \eqref{AP4} we have 
\begin{equation*}
\|\bar F\|_{L_p(\Omega)} \leq C[ \|F\|_{L_p(\Omega)}+\|\de_{x_1}u\|_{L_p(\Omega)}+\|u,B\|_{W^{1-1/p}_p(\Gamma)}],
\end{equation*}
therefore denoting  
\begin{equation}
-(2\mu+\nu)\div u + Kw=H \label{AP4b}
\end{equation}
we have for any $\delta>0$
\begin{align} \label{AP5}
\|H\|_{W^s_p(\Omega)} \leq &\delta \|\nabla H\|_{L_p(\Omega)}+C(\delta)\|H\|_{L_2(\Omega)}\leq
\delta [\|F\|_{L_p(\Omega)}+\|\de_{x_1}u\|_{L_p(\Omega)}+\|u,B\|_{W^{1-1/p}_p(\Gamma)}]\nonumber\\
&+C(\delta)\|\nabla u,w\|_{L_2(\Omega)}.
\end{align}

%Looking at the potential part of $(\ref{system_lin})_1$ we get
%\begin{equation}\label{AP4}
% \nabla (-(2\mu+\nu) \div u + K w) = P_{\nabla} (F-\de_{x_1}u),
%\end{equation}
%where $P_{\nabla} F$ is the potential part of $F$. 
%Therefore denoting
%\begin{equation}\label{AP5}
% -(2\mu+\nu) \div u + K w=F_2
%\end{equation}
%fulfilling the following inequ
%\begin{equation}\label{AP6}
% \|f_3\|_{W^s_p(\Omega)} \leq C(\|F\|_{L_p(\Omega)} + \|\nabla u,w\|_{L_2(\Omega)}).
%\end{equation}
%The last terms links possibility of removing the gradients (\ref{AP4}) with the energy norm in order to compute the average of the l.h.s. of (\ref{AP5}).
%{\bf describe more precisely...}

\bigskip

{\bf Step III}.
Adding $(\ref{system_lin})_2$ to (\ref{AP4b}) with suitable scaling we obtain
\begin{equation} \label{AP7}
\d_{x_1}w + U\cdot \nabla w + \bar K w = \bar G
\end{equation}
with $\bar K = \frac{K}{2\mu+\nu}$ and 
$$
\|\bar G\|_{W^s_p(\Omega)} \leq \|G\|_{W^s_p(\Omega)}+ [RHS \; of \; \eqref{AP5}].
$$
Here we meet the main mathematical challange of our result, we have to solve this transport like problem in a domain which touches the characteristics at the ends
of  $\Gamma_{in}$. Its solvability is given by Proposition \ref{lem_trans}. 
In particular, the estimate \eqref{est_trans} yields
\begin{equation} \label{AP8}
 \|w\|_{W^s_p(\Omega)} \leq C (\|G\|_{W^s_{p}(\Omega)} 
 + \|w_{in}\|_{W^s_p(\Gamma_{in}) \cap W^r_p(\Gamma_s)}).
\end{equation}
\begin{rem} \label{rem_indep2}
Notice that Remark \ref{rem_indep} implies that we can assume that the constant in \eqref{AP8} 
does not depend on $\|U\|_{W^1_\infty}$.  
\end{rem}
{\bf Step IV.} We collect the elements of our estimation. 
In order to get the information about the velocity we use estimates 
for $\rot u$ given by (\ref{AP3}) and 
about $\div u$ coming from \eqref{AP4b},\eqref{AP5} and \eqref{AP8}. We obtain
\begin{align} \label{AP9}
& \|u\|_{W^{1+s}_p(\Omega)} \leq C(\|\rot u\|_{W^s_p(\Omega)}+\|\div u\|_{W^s_p(\Omega)}) \nonumber \\
&\leq C( \|F\|_{L_p(\Omega)} + \|G\|_{W^s_p(\Omega)} + \|B\|_{W^{1-1/p}(\Gamma)} + \|\nabla u, w\|_{L_2(\Omega)}
&+ \|\d_{x_1}u \|_{L_p} +\|u\|_{W^{1-1/p}(\Gamma)} )\nonumber \\
&+\delta (\|u\|_{W^{1+s}_p(\Omega)}+\|w\|_{W^s_p(\Omega)}).
\end{align}
The $L_2$ term on the rhs is treated with the energy estimate \eqref{ene}. 
To the boundary term $\|u\|_{W^{1-1/p}(\Gamma)}$ 
we apply first trace theorem and then interpolation inequality \eqref{int} to get
\begin{equation} \label{AP10}
\|u\|_{W^{1-1/p}_p(\Gamma)} \leq C \|u\|_{W^s_p(\Omega)} \leq \delta \|u\|_{W^{1+s}_p(\Omega)} + C(\delta) \|\nabla u\|_{L_2(\Omega)}.
\end{equation}
The term $\d_{x_1}u$ is treated similarly with interpolation inequality.
We conclude \eqref{est_lin} and complete the proof. 
\qed
Now it is a matter of standard theory to show the existence for the linear system \eqref{system_lin}. 
\begin{prop} \label{prop:lin}
Let the rhs of \eqref{system_lin} be of regularity specified in \eqref{reg_rhs}
with $s,p$ and the boundary data satisfying the assumptions of Theorem 1. Then
there exists a unique solution to \eqref{system_lin} satisfying the estimate 
\eqref{est_lin}.
\end{prop}
\emph{Proof.} As we have the a priori estimate, the proof of existence is standard. 
We start with showing the weak solutions using the Galerkin method and the energy estimate 
\eqref{ene}. Then, using our {\emph a priori} estimate we show that our solution has the 
required regularity. \qed

\subsection{Estimate for the nonlinear system}
We are now ready to close the estimate in $W^{1+s}_p \times W^s_p$ for the nonlinear system. 
To this end we combine the linear estimate
\eqref{est_lin} with the bounds \eqref{a7} on the rhs of the nonlinear problem obtaining
%We look at the rhs of (\ref{AP8}), take the first term
%\begin{equation}\label{AP9}
% \|F_1\|_{L_p(\Omega)} \leq C(\|\d_{x_1} u \|_{L_p(\Omega)} + \|w\|_{L_\infty(\Omega)}\|\nabla u\|_{L_p(\Omega)} + E(\|\nabla u,u,w\|_{L_p(\Omega)} +1)),
%\end{equation}
%where $E$ is small in terms of the data.
%Here the most non-obvious term is the first one, we apply the interpolation inequality 
%\eqref{int} getting (for some $\theta \in (0,1)$)
%\begin{equation}\label{AP10}
% \|\d_{x_1} u \|_{L_p(\Omega)} \leq C\|u\|_{W^{1+s}_p(\Omega)}^\theta\|\nabla u\|_{L_2(\Omega)}^{1-\theta}\leq \epsilon \|u\|_{W^{1+s}_p(\Omega)} + C\|\nabla u \|_{L_2(\Omega)}.
%\end{equation}
%Next, we find
%\begin{equation}\label{AP11}
% \|f_2\|_{W^s_p(\Omega)}\leq C(\|w\|_{W^s_p(\Omega)}^2 + E(\|w\|_{W^s_p(\Omega)} +1))
%\end{equation}
%and 
%\begin{equation}\label{AP12}
% \|G\|_{W^s_p(\Omega)}\leq C(\|w\|_{W^s_p(\Omega)}\|u\|_{W^{1+s}_p(\Omega)} + E(\|w\|_{W^s_p(\Omega)} +1)).
%\end{equation}
%
\begin{equation*}\label{AP15}
\|u\|_{W^{1+s}_p(\Omega)} + \|w\|_{W^s_p(\Omega)} \leq
E( \|u\|_{W^{1+s}_p(\Omega)} + \|w\|_{W^s_p(\Omega)}+1)
+C( \|u\|_{W^{1+s}_p(\Omega)} + \|w\|_{W^s_p(\Omega)})^2,
\end{equation*}
where $E$ is a small constant dependent on the data and we can assume (see Remark \ref{rem_indep2})
that $C$ does not depend on $\|u\|_{W^1_\infty(\Omega)}$. For sufficiently small data we conclude
\begin{equation} \label{AP16}
\|u\|_{W^{1+s}_p(\Omega)} + \|w\|_{W^s_p(\Omega)} \leq 
C( \|u\|_{W^{1+s}_p(\Omega)} + \|w\|_{W^s_p(\Omega)}^2 + E,
\end{equation}
where $C$ does not depend on $(\|w\|_{W^s_p(\Omega)},\|u\|_{W^{1+s}_p(\Omega)})$.
This estimate will be crucial in showing the existence and uniqueness of solutions in the next section.
\begin{rem}
Notice that the estimate \eqref{AP16} holds without any assumptions on relation between $\mu$ and $K$. 
However, this relation will come into play in the following section.
\end{rem}

\section{Proof of Theorem 1}
In order to prove our main result we apply an iterative scheme.
The idea is to combine the estimate \eqref{AP16} with the Cauchy condition in some weaker space. 
In \cite{P1},\cite{PP} this space was $H^1 \times L_2$. However, then we need some information on the gradient of the density
which is not available in our framework. Here we overcome this obstacle showing the convergence 
in $L_2 \times H^{-1}$. However, then we have to estimate $\|\pi(\rho^n)-\pi(\rho^{n-1})\|_{H^{-1}}$ in terms 
of $\|\rho^n-\rho^{n-1}\|_{H^{-1}}$ and for this purpose we have to assume \eqref{plin}. 
A similar approach was applied in the context of uniqueness of weak solutions in the nonstationary 
case in \cite{Hoff} where the constraint \eqref{plin} also appeared.
Under this assumption, denoting
$$
v^{n+1}=u^{n+1}+u_0+\bar v
$$
we can define the sequence of approximations in the following way  
\begin{equation} \label{uniq:00}
\begin{array}{lcr}
\div(w^{n+1} v^{n+1} \otimes v^{n+1})-\mu \Delta u^{n+1}-(\mu+\nu)\nabla \div u^{n+1}+K \nabla w^{n+1}=0 
& \mbox{in} & \Omega,\\[3pt]
w^{n+1}_{x_1}+\div((u^n+u_0)w^{n+1})=-\div (u^n+u_0)
& \mbox{in}& \Omega,\\[3pt]
n\cdot 2\mu {\bf D}(u^{n+1})\cdot \tau +f u^{n+1} \cdot \tau = B
&\mbox{on} & \Gamma, \\[3pt]
n\cdot u^{n+1} = 0 & \mbox{on} & \Gamma,\\[3pt]
w^{n+1}=w_{in} & \mbox{on} & \Gamma_{in}.
\end{array}
\end{equation}
We start with the following auxiliary result
\begin{lem}  
Assume that $\phi$ solves 
\begin{equation} \label{uniq:aux}
%\begin{array}{lr}
(\de_{x_1}+U\cdot\nabla)\phi=f, \quad \phi|_{\Gamma_{out}}=0
%\end{array}
\end{equation}
with $f \in H^1_0(\Omega)$ and sufficiently small $U \in W^1_\infty(\Omega)$ such that 
\begin{equation} \label{uniq:0}
U\cdot n|_{\Gamma_{in}}<0, \quad U\cdot n|_{\Gamma_{out}}>0. 
\end{equation} 
Then 
\begin{equation} \label{uniq:0b}
\|\nabla \phi\|_{L_2(\Omega)}\leq C_\Omega \|f\|_{H^1_0(\Omega)}
\end{equation}
where  $C_\Omega \sim {\rm max}\{|\Omega|,|\Omega|^2\}$.
\end{lem}
\emph{Proof.} 
The boundary condition \eqref{uniq:aux}$_2$ implies 
\begin{equation} \label{uniq:2}
0=\de_{\tau}\phi|_{\Gamma_{out}}=\tau_1\phi_{x_1}+\tau_2\phi_{x_2}|_{\Gamma_{out}}.
\end{equation}
On the other hand we have
\begin{equation} \label{uniq:3} 
(1+u^1)\phi_{x_1}+u^2\phi_{x_2}|_{\Gamma_{out}}=0.
\end{equation}
Subtracting \eqref{uniq:2} multiplied by $u^2$ from \eqref{uniq:3} multiplied by $\tau^2$ w get 
$$
0=[(1+u^1)\tau^2-u^2\tau^1)]\phi_{x_1}=-[(1+u^1)n^1+u^2n^2]\phi_{x_1}=0,
$$   
therefore $\phi_{x_1}|_{\Gamma_{out}}=0$ due to \eqref{uniq:0}. Combining this identity 
with \eqref{uniq:2} we conclude 
\begin{equation} \label{uniq:4}
\nabla \phi|_{\Gamma_{out}}=0.
\end{equation}
Next we differentiate \eqref{uniq:aux}$_1$ wrt $x_2$ and multiply by $\phi_{x_2}$.
Denoting $V=[1+U^1,U^2]$ we obtain
\begin{equation} \label{uniq:5}
\frac{1}{2}V \cdot \nabla \phi_{x_2}^2=f_{x_2}\phi_{x_2}-V_{x_2}\phi_{x_2}\cdot\nabla\phi.
\end{equation} 
Now let us denote 
$$
\Omega_{a}=\Omega \cap \{x_1>a\}, \quad I_a = \Omega \cap \{x_1=a\}.
$$
Integrating \eqref{uniq:5} over $\Omega_{a}$ we get 
\begin{equation} \label{uniq:6}
\frac{1}{2}\int_{\de \Omega_{x_1}} (V\cdot n) \phi_{x_2}^2\,d\sigma=
\int_{\Omega_{x_1}}\Big[ \frac{1}{2}\phi_{x_2}^2\div V 
+f_{x_2}\phi_{x_2} + V_{x_2}\phi_{x_2}\cdot\nabla\phi \Big] \,dx.
\end{equation} 
Notice that by \eqref{uniq:4} and the condition $V \cdot n|_{\Gamma_{in}}<0$ we have 
$$
\int_{\de \Omega \, \cap \, \de \Omega_{x_1}} (V\cdot n) \phi_{x_2}^2\,d\sigma \leq 0.
$$
Therefore the smallness of $U^1$ implies 
$$
\int_{\de \Omega_{x_1}} (V\cdot n) \phi_{x_2}^2\,d\sigma =
-\int_{I_{x_1}}(1+U^1)\phi_{x_2}^2\,dx_2 + \int_{\de \Omega \, \cap \, \de \Omega_{x_1}}(V\cdot n) \phi_{x_2}^2\,d\sigma \leq 
- C_U \int_{I_{x_1}} \phi_{x_2}^2 \,d\sigma
$$
for some $C_U>0$. The latter inequality combined with \eqref{uniq:6} gives for any $\ep>0$ 
$$
{\rm sup}_{x_1} \int_{I_{x_1}}\phi_{x_2}^2\,dx_2 \leq \| \nabla U \|_{L_\infty(\Omega)} \|\nabla \phi\|_{L_\infty(\Omega)}
+\ep \|\phi_{x_2}\|_{L_2(\Omega)}^2 + \frac{C}{\ep} \|\nabla f\|_{L_2(\Omega)}^2.
$$
Using again smallness of $U$ we obtain  
\begin{equation} \label{uniq:7}
\|\phi_{x_2}\|_{L_2(\Omega)}^2 \leq |\Omega| [(\ep+\|\nabla U\|_{L_\infty(\Omega)}) \|\nabla \phi\|_{L_2(\Omega)}^2+C(\ep)\|\nabla f\|_{L_2(\Omega)}^2 ].
\end{equation}
Finally from \eqref{uniq:aux}$_1$ we have 
\begin{equation*} 
\|\phi_{x_1}\|_{L_2(\Omega)}^2 \leq \|\frac{f}{1+U^1}\|_{L_2(\Omega)}^2+\| \frac{U^2}{1+U^1} \|_{L_\infty(\Omega)}\|\phi_{x_2}\|_{L_2(\Omega)}^2. 
\end{equation*}
Combining this inequality with \eqref{uniq:7} we get  
\begin{equation*}
\|\nabla \phi\|_{L_2(\Omega)}^2\leq [|\Omega|(\ep+\nabla U\|_{L_\infty(\Omega)})+C\|U\|_{L_\infty(\Omega)}]\|\nabla \phi\|_{L_2(\Omega)}^2
+[\frac{C|\Omega|}{\ep}+C_p]\|\nabla f\|_{L_2(\Omega)}^2
\end{equation*}
where $C_p$ is the constant from the Poincar\'e inequality.
Now, provided $\|U\|_{W^1_\infty}$ is sufficiently small, to close the estimate we need $\ep \sim \frac{1}{|\Omega|}$.
Taking into account that $C_p \sim |\Omega|$ we conclude \eqref{uniq:0b}. 
\qed  

\noindent
The following proposition implies convergence of the sequence $(u^n,w^n)$ in $L_2 \times H^{-1}$.
\begin{prop} \label{uniq:prop}
Assume that $\mu$ is sufficiently large compared to $|\Omega|$, $\|\bar v\|_{L_\infty}$
and $K$.
Then the sequence $(u^n,w^n)$ defined by \eqref{uniq:00} satisfies 
\begin{equation} \label{cc}
\|u^{n+1}-u^n\|_{L_2(\Omega)}+\|w^{n+1}-w^n\|_{H^{-1}(\Omega)}\leq M[\|u^n-u^{n-1}\|_{L_2(\Omega)}+\|w^n-w^{n-1}\|_{H^{-1}(\Omega)}].
\end{equation}
with $M<1$.
\end{prop}
\begin{rem} In order to track the dependence of $\mu$ on $\|\bar v\|_{L_\infty}$ we consider 
again \eqref{const} with a general constant $v^*$.
\end{rem} 
\emph{Proof.}  
Subtracting the equations \eqref{uniq:00}$_2$ for two consecutive steps we get 
\begin{equation*}
(w^{n+1}-w^n)_{x_1}+\div [(u^n+u_0)(w^{n+1}-w^n)]=-\div( (w^n+1)(u^n-u^{n-1})). 
\end{equation*} 
We test this equation with $\phi$ given by \eqref{uniq:aux} with $U=u^n+u_0$ and $f$ such that 
\begin{equation} \label{uniq:10}
\Delta f = w^{n+1}-w^n.
\end{equation}
We obtain 
\begin{equation*}
\int_{\Omega}(u^n-u^{n-1})(w_n+1)\cdot \nabla \phi\dx = 
-\int_{\Omega}(w^{n+1}-w^n)(\de_{x_1}+(u^n+u_0)\cdot\nabla)\phi\dx=-\int_{\Omega}f \Delta f\dx=\|w^{n+1}-w^n\|_{H^{-1}}^2.
\end{equation*}
Therefore by \eqref{uniq:0b} we obtain 
\begin{equation*}
\|w^{n+1}-w^n\|_{H^{-1}}^2 \leq \|w^n+1\|_{L_\infty}\|u^n-u^{n-1}\|_{L_2}\|\nabla \phi\|_{L_2} \leq 
C_1 C_{\Omega} \|u^n-u^{n-1}\|_{L_2} \|\nabla f\|_{L_2}   
\end{equation*}
where $C_{\Omega}$ is the constant from \eqref{uniq:0b}.
Now, since $\|\nabla f\|_{L_2} \leq C\|w^{n+1}-w^n\|_{H^{-1}}$, we conclude 
\begin{equation} \label{uniq:11}
\|w^{n+1}-w^n\|_{H^{-1}(\Omega)} \leq C_2 C_{\Omega} \|u^n-u^{n-1}\|_{L_2(\Omega)}. 
\end{equation}
Now we have to estimate $\|u^n-u^{n-1}\|_{L_2(\Omega)}$ in terms of $\|w^n-w^{n-1}\|_{H^{-1}(\Omega)}$.  
Subtracting the equations \eqref{uniq:00}$_1$ for two consecutive steps we obtain 
\begin{align} \label{uniq:12}
&-\mu \Delta (u^{n+1}-u^n)-(\mu+\nu)\nabla \div (u^{n+1}-u^n)= -K \nabla (w^{n+1}-w^n) - \div [(w^{n+1}-w^n)(v^{n+1}\otimes v^{n+1})] \nonumber\\
&-\div [ w^n( v^{n+1} \otimes (v^{n+1}-v^n)+(v^{n+1}-v^n)\otimes v^n ) ].
\end{align}
Notice that as $u^{n+1}-u^n$ satisfy homogeneous boundary conditions we have in the weak sense   
$$
\int_{\Omega} \psi \cdot [\mu \Delta u + (\mu+\nu)\nabla \div u]dx=
\int_{\Omega} u \cdot [\mu \Delta \psi + (\mu+\nu)\nabla \div \psi]dx
$$
for any $\psi \in H^2: \psi|_{\Gamma}=0$.
Therefore testing \eqref{uniq:12} with $\psi$ being a solution to the problem 
\begin{equation*}
-\mu \Delta \psi - (\mu+\nu)\nabla \div \psi = \mu (u^{n+1}-u^n), \quad \psi|_{\Gamma}=0
\end{equation*} 
we get
\begin{align*}
&\mu \|u^{n+1}-u^n\|_{L_2(\Omega)}^2 = \int_{\Omega}(w^{n+1}-w^n)(K \div \psi+v^{n+1}\otimes v^{n+1}:\nabla \psi) dx\nonumber\\
&+\int_{\Omega} w^n[ v^{n+1}\otimes(v^{n+1}-v^n)+(v^{n+1}-v^n)\otimes v^n]:\nabla\psi dx \leq \nonumber\\[3pt]
&\leq C_v[\|w^{n+1}-w^n\|_{H^{-1}(\Omega)}\|\nabla \psi\|_{L_2(\Omega)} + \|u^{n+1}-u^n\|_{L_2(\Omega)}\|\nabla \psi\|_{L_2(\Omega)}\leq\nonumber\\[3pt]
&\leq C_{v,K}[ \|w^{n+1}-w^n\|_{H^{-1}(\Omega)}\|u^{n+1}-u^n\|_{L_2(\Omega)}+\|u^{n+1}-u^n\|_{L_2(\Omega)}^2 ],  
\end{align*} 
where $C_{v,K}=C_{v,K}(\|v^{n+1}\|_{L_\infty(\Omega)},K)$. In the above passages we have used the facts that 
$v^{n+1}-v^n=u^{n+1}-u^n$ and $\|\nabla \psi\|_{L_2(\Omega)} \leq C \|u^{n+1}-u^n\|_{L_2(\Omega)}$.
Therefore, assuming $\mu>C_v$ we obtain
\begin{equation*}
\|u^{n+1}-u^n\|_{L_2(\Omega)}\leq \frac{C_{v,K}}{\mu-C_{v,K}}\|w^{n+1}-w^n\|_{H^{-1}(\Omega)}.
\end{equation*}
Combining this estimate with \eqref{uniq:11} we conclude \eqref{cc} 
with $M=\frac{C_2 C_\Omega C_{v,K}}{\mu-C_{v,K}}$. Therefore $M<1$
provided the viscosity is sufficiently large compared to $|\Omega|$,$\|\bar v\|_{L_\infty(\Omega)}$ and $K$. 
\qed
The solvability of the linear system established in Section 3 implies that the sequence 
\eqref{uniq:00} is well defined. Moreover, denoting 
$$
A_n = \|u^n\|_{W^{1+s}_p(\Omega)}+\|w^n\|_{W^s_p(\Omega)},
$$
by \eqref{AP16} we have
$$
A_{n+1} \leq C A_n^2 + E 
$$
and we can assume $E \leq \frac{1}{2C}$. Then the sequence starting with 
$(u^0,w^0)=([0,0],0)$ satisfies 
\begin{equation} \label{uniq:10}
A_n \leq 2 E
\end{equation}
where $E$ is the constant from \eqref{AP16}. Proposition \ref{uniq:prop} implies that 
$$
(u^n,w^n) \to (u,w) \quad {\rm in} \quad L_2(\Omega) \times H^{-1}(\Omega). 
$$
On the other hand, \eqref{uniq:10} yields up to a subsequence 
$(u^n,w^n) \rightharpoonup (\bar u, \bar w)$ in $W^{1+s}_p(\Omega) \times W^s_p(\Omega)$. 
By the definition of $(u^n,w^n)$, the limit $(u,w)$ is a solution to \eqref{system}-\eqref{FG} 
in the sense of Definition \ref{def_sol1}. The uniqueness is shown in the same way as Proposition \ref{uniq:prop}.
Namely, taking two solutions $(u^1,w^1)$ and $(u^2,w^2)$ for the same data we show
$$
\|u^1-u^2\|_{L^2(\Omega)}+\|w^1-w^2\|_{H^{-1}(\Omega)}=0.
$$   
\qed

\section{Concluding remarks}
The solutions considered here are located somehow between weak and "traditional" regular 
solutions satisfying the equations almost everywhere. The result for the steady 
transport equation given by Lemma \ref{lem_trans} is obtained for a general class of boundary singularities 
showing that the choice of the fractional Sobolev-Slobodetskii spaces is in a sense 
natural for the problem under consideration. The price we pay 
is that we have to assume linearity of the pressure. Getting rid of this constraint 
seems an interesting open problem. It is likely that the existence itself could 
be shown for more general pressure laws using for example approximation with more regular solutions which give some information 
about the gradient of the density. However, such result would be highly technical and not really meaningful 
without uniqueness which is more challenging and seem to require some novel approach to treat more general pressure laws.   
Finally we shall mention that the assumed $C^3$ regularity of the domain required to solve 
an elliptic problem appearing in the estimate for the velocity is not optimal and could be relaxed at a price of additional technicalities. 
However, as we are rather interested in a careful investigation of the boundary singularity in the stationary transport 
equation which is independent of global regularity of the boundary (for example a $C^\infty$ boundary 
can fail to satisfy \eqref{flat}), we keep this regularity assumption.

\smallskip
  
\noindent
{\bf Acknowledgements.} This work has been supported by Ideas Plus grant ID 2011 0006 61. The authors would like 
to thank the anonymous Referee for a careful lecture of the manuscript and numerous remarks which contributed 
to improving the clarity of the proofs.

\end{document}